\documentclass[12pt]{amsart}
\usepackage{amscd}
\usepackage{amssymb}
\usepackage{amsmath}
\usepackage{amsfonts}

\newtheorem{lemme}{Lemme}[section]
\newtheorem{proposition}[lemme]{Proposition}
\newtheorem{thr}[lemme]{Th\'eor\`eme}
\newtheorem{corollaire}[lemme]{Corollaire}

\newcommand\dem{\noindent{\it D\'emonstration.}\ }
\newcommand\findem{\hfill$\square$ \bigskip}

\newcommand\FF{{\mathbb F}}
\newcommand\LL{{\mathbb L}}
\newcommand\PP{{\mathbb P}}
\newcommand\QQ{{\mathbb Q}}
\newcommand\ZZ{{\mathbb Z}}

\newcommand\Gr{{\mathcal Q}}
\newcommand\Lie{\mathrm{Lie}}
\newcommand\Hom{\mathrm{Hom}}
\newcommand\Gm{{\mathbb G}_m}

\newcommand\calO{{\mathcal O}}
\newcommand\calA{{\mathcal A}}
\newcommand\calC{{\mathcal C}}
\newcommand\calH{{\mathcal H}}
\newcommand\calL{{\mathcal L}}

\newcommand\rmH{{\rm H}}
\newcommand\rmR{{\rm R}}
\newcommand\Ql{{\bar{\mathbb Q}}_\ell}
\newcommand\vp{\varpi}
\newcommand\rta{\rightarrow}
\newcommand\Proj{{\rm Proj\,}}
\newcommand\smv{{\scriptscriptstyle\vee}}
\newcommand\tildetimes{\,{\tilde\times}\,}

\newcommand\cf{\mbox{\it c.f.\, }}
\newcommand\loc{\mbox{\it loc.\,cit.}}
\newcommand\deuxind[2]{{\buildrel{\scriptstyle #1}\over{#2}}}

\newcommand\hfld[2]{\smash{\mathop{\hbox to 12mm{\rightarrowfill}}
\limits^{\scriptstyle#1}_{\scriptstyle#2}}}
\newcommand\hflg[2]{\smash{\mathop{\hbox to 12mm{\leftarrowfill}}
\limits^{\scriptstyle#1}_{\scriptstyle#2}}}

\begin{document}
\title[Formule de Casselman-Shalika g\'eom\'etrique]
{R\'esolutions de Demazure affines et 
formule de Casselman-Shalika g\'eom\'etrique}
\author{B.C. Ng\^o et P. Polo}
\date{}

\begin{abstract}
We prove a conjecture of Frenkel, Gaitsgory, Kazhdan, Vilonen,
related to Fourier coefficients of spherical perverse sheaves on
the affine Grassmannian associated to a split reductive group.
\end{abstract}

\maketitle

\section*{Introduction}

Soit $G$ un groupe r\'eductif connexe d\'eploy\'e sur un corps
fini $k={\FF}_q$.
Notons $\calO=k[[\varpi]]$ l'anneau des
s\'eries formelles \`a coefficients dans $k$ et $F$ le corps des fractions
de $\calO$.
Soit $K=G(\calO)$ le sous-groupe compact maximal standard de $G(F)$.
Pour tout cocaract\`ere dominant $\lambda$ de $G$, il est possible de construire
un $k$-sch\'ema
projectif ${\bar\Gr}_\lambda$ dont l'ensemble des $k$-points est
$$
{\bar\Gr}_\lambda(k)=\coprod_{\lambda'\leq\lambda}K\varpi^{\lambda'}K/K ,
$$
sur lequel agit le groupe $K$,
vu comme un $k$-groupe alg\'ebrique de dimension infinie, \`a travers un
quotient de type fini.
Cette action induit une stratification en orbites
${\bar\Gr}_\lambda=\coprod_{\lambda'\leq\lambda}\Gr_{\lambda'}$
parmi lesquelles $\Gr_\lambda$ est l'orbite ouverte.
De plus, les ${\bar\Gr}_\lambda$ s'organisent en une famille inductive
dont la limite
est le r\'eduit associ\'e \`a la Grassmannienne affine $\Gr$ d\'efinie
comme dans \cite{Beauville-Laszlo,Laszlo-Sorger}.

Le sch\'ema ${\bar\Gr}_\lambda$ n'\'etant pas lisse en g\'en\'eral,
pour un nombre premier $\ell$ diff\'erent de la caract\'eristique de $k$,
il est naturel
de consid\'erer le complexe d'intersection $\ell$-adique
$$
\calA_\lambda={\rm IC}({\bar\Gr}_\lambda,\Ql) ,
$$
qui est $K$-\'equivariant.
La fonction trace de Frobenius associ\'ee \`a ce faisceau pervers :
$$
A_\lambda(x)={\rm Tr}({\rm Fr}_q,(\calA_\lambda){}_x) ,
$$
d\'efinie sur l'ensemble des $k$-points
de ${\bar\Gr}_\lambda$, peut \^etre vue comme
un \'el\'ement de l'alg\`ebre de Hecke non ramifi\'ee $\calH$ des fonctions
\`a valeurs dans $\Ql$,
\`a support compact dans
$G(F)$, qui sont invariantes \`a gauche et \`a droite par $K$. Lorsque
$\lambda$ parcourt
le c\^one des cocaract\`eres dominants,
ces fonctions $A_\lambda$ forment une base de $\calH$.

Soit $G^\smv$ le groupe d\'efini sur ${\Ql}$ dont la donn\'ee radicielle
est duale de celle de $G$.
Dans \cite{Satake}, Satake a construit un isomorphisme
canonique entre l'alg\`ebre $\calH$ et
l'alg\`ebre des fonctions r\'eguli\`eres
sur $G^\smv$ qui sont ${\rm Ad}(G^\smv)$-invariantes. D'apr\`es Lusztig et
Kato, voir \cite{Lusztig}, \cite{Kato}, la transformation
de Satake de $A_\lambda$ est \'egale, \`a un signe pr\`es, au caract\`ere de
$V(\lambda)$, la repr\'esentation irr\'eductible de plus haut poids
$\lambda$ de $G^\smv$. Plus r\'ecemment, Ginzburg \cite{Ginzburg} et
Mirkovic, Vilonen \cite{Mirkovic-Vilonen} ont mis en lumi\`ere
une \'equivalence tannakienne entre la cat\'egorie des faisceaux pervers
$K$-\'equi\-va\-riants semi-simples
sur $\Gr$ munie d'un produit de convolution et la cat\'egorie des
repr\'esen\-ta\-tions alg\'ebriques de $G^\smv$ munie du produit
tensoriel. Le th\'eor\`eme de Lusztig-Kato serait le reflet au niveau des
objets simples de cette \'equivalence, via la formule
des traces de Grothendieck
\cite{Grothendieck}.

Les termes constants ainsi que les coefficients de Fourier des fonctions
$A_\lambda$ sont remarquablement simples. Soit
$B = TU$ un sous-groupe de Borel de $G$
et $\rho$ la demi-somme des racines de $T$ dans $\Lie(U)$.
D'apr\`es Lusztig et Kato, l'int\'egrale terme constant est \'egale \`a
$$
\int_{U(F)}A_\lambda(x\varpi^\nu)\, dx=(-1)^{2\langle\rho,\nu\rangle}
q^{\langle\rho,\nu\rangle}m_\lambda(\nu) ,
$$
o\`u $m_\lambda(\nu)$ est
la dimension de l'espace de poids $\nu$ dans $V(\lambda)$.
Parall\`element, pour $\theta:U(F)\rightarrow \QQ_\ell^\times$ un
carac\-t\`ere g\'en\'erique,
de conducteur $U(\calO)$, Frenkel, Gaitsgory, Kazhdan et Vilonen ont
d\'emontr\'e dans \cite{FGKV}, que
$$\int_{U(F)}A_\lambda(x\varpi^\nu)\theta(x)\, dx=0$$
si $\nu\not=\lambda$ et
$$
\int_{U(F)}A_\nu(x\varpi^\nu)\theta(x)\, dx=(-1)^{2\langle\rho,\nu\rangle}
q^{\langle\rho,\nu\rangle}.
$$
Leur d\'emonstration s'appuie sur la formule explicite de Cassel\-man-Shalika
des valeurs des fonctions de Whittaker non
ramifi\'ees \cite{Shintani}, \cite{Casselman-Shalika}.

L'objet principal de notre travail est de d\'emontrer les \'enonc\'es
g\'eo\-m\'e\-tri\-ques
sous-jacents \`a ces r\'esultats. Pour tout cocaract\`ere $\nu$, il est
possible de d\'efinir un sous-ind-sch\'ema
$S_\nu\subset\Gr$ tel que $$S_\nu(k)=U(F)\vp^\nu G(\calO)/G(\calO).$$ Il
s'agit de d\'emontrer que le complexe
$$
{\rm R}\Gamma_c(S_\nu\otimes_k{\bar k},{\mathcal A}_\lambda)
$$
est concentr\'e en degr\'e $2\langle\rho,\nu\rangle$ et que
l'endomorphisme de Frobenius agit dans son $\rmH^{2\langle\rho,\nu\rangle}$
comme la multiplication par
$q^{\langle\rho,\nu\rangle}$. Cet
\'enonc\'e est d\^u \`a Mirkovic et Vilonen lorsque le corps de base est
${\mathbb C}$ et joue un r\^ole important dans
l'\'equivalence tannakienne mentionn\'ee plus haut. Il peut aussi \^etre
consid\'er\'e comme une interpr\'etation g\'eom\'etrique
partielle du th\'eor\`eme de Lusztig-Kato.

Lorsque $\nu$ est dominant, on peut d\'efinir un morphisme
$h:S_\nu\rta{\mathbb G}_a$
tel que $\theta(x)=\psi(h(x))$, o\`u $\psi:k\rta\Ql^\times$ est un
caract\`ere additif non trivial de $k$.
On d\'emontre que le
complexe
$$
{\rm R}\Gamma_c(S_\nu\otimes_k{\bar k},{\mathcal A}_\lambda\otimes
h^*{\mathcal L}_\psi)
$$
est nul si $\nu\not=\lambda$. Dans le cas $\nu=\lambda$,
il est isomorphe \`a $\Ql(-\langle\rho,\nu\rangle)$ plac\'e en degr\'e
$2\langle\rho,\nu\rangle$.
Cet \'enonc\'e \'etait une conjecture de
Frenkel, Gaitsgory, Kazhdan et Vilonen
\cite{FGKV}. Comme expliqu\'e dans \loc,
il fournit une d\'emonstration g\'eom\'etrique du th\'eor\`eme de
Casselman-Shalika.
Il pourrait \'egalement fournir des
r\'esultats pour les groupes tordus.

Voici l'organisation de l'article. Apr\`es avoir rappel\'e, dans la section 2,
des r\'esultats connus sur la Grassmannienne affine,
nous \'enon\c cons les r\'esultats principaux (th\'eor\`emes 3.1 et 3.2)
dans la section 3.

La d\'emonstration de ces th\'eor\`emes occupe
le reste de l'article.
Elle repose
sur l'\'etude de la g\'eom\'etrie de certaines r\'esolutions, form\'ees \`a
partir des vari\'et\'es ${\bar\Gr}_\lambda$ les plus simples,
qui correspondent au cas o\`u $\lambda$ est minuscule ou
quasi-minuscule. Cette strat\'egie est d\'ej\`a
utilis\'ee dans \cite{Ngo}, o\`u la conjecture de
\cite{FGKV} a \'et\'e d\'emontr\'ee pour le groupe ${\rm GL}_n$.

De fa\c con plus d\'etaill\'ee,
dans les sections 4 et 5, nous d\'emontrons de fa\c con g\'eom\'etrique
deux \'enonc\'es (lemmes 4.2 et 5.2)
concernant les intersections $S_\nu\cap{\bar\Gr}_\lambda$, qui sont
probablemment bien
connus, mais que nous n'avons su trouver, sous cette forme, dans la
litt\'erature.
Le lemme 5.2 nous permet de d\'emontrer les th\'eor\`emes 3.1 et 3.2
dans le cas o\`u $\nu$ et
$\lambda$ sont conjugu\'es par un \'el\'ement du groupe de Weyl.

Nous signalons, au passage, un \'enonc\'e (proposition 4.3) sur les
ca\-rac\-t\'e\-ris\-tiques d'Euler-Poincar\'e
$\chi_c(S_\nu\cap\Gr_\lambda)$ qui peut \^etre consid\'er\'e comme une
inter\-pr\'e\-ta\-tion g\'eom\'etrique d'un r\'esultat de Lusztig
\cite[6.1]{Lusztig}.

Nous \'etudions ensuite,
dans les sections 6--8, la g\'eom\'etrie des vari\'et\'es
${\bar\Gr}_\lambda$ dans les cas les plus simples, c.\`a.d.
lorsque $\lambda$ est minuscule (section 6), ou quasi-minuscule (sections 7
et 8).

Si $\lambda$ est minuscule, alors ${\bar\Gr}_\lambda$ est
\'egal \`a
$\Gr_\lambda$ et est isomorphe au sch\'ema $G/P$ des sous-groupes de $G$
conjugu\'es \`a
un certain sous-groupe parabolique $P$; de plus, seuls les $\nu$
conjugu\'es \`a $\lambda$ interviennent, si
bien que les \'enonc\'es 3.1 et 3.2 d\'ecoulent dans ce cas du lemme 5.2.

Si $\lambda$ est quasi-minuscule, alors le point base de $\Gr$ est l'unique
point singulier de ${\bar\Gr}_\lambda$.
L'ouvert compl\'ementaire de ce point, l'orbite $\Gr_\lambda$, est un
fibr\'e en droites
au-dessus d'un $G/P$. Nous construisons une r\'esolution (lemme 7.3) de
${\bar\Gr}_\lambda$ par un fibr\'e en droites projectives
au-dessus de $G/P$, qui nous permettra de d\'emontrer, dans la section 8,
les th\'eor\`emes 3.1 et 3.2 dans ce cas.

Dans la section 9, on consid\`ere certaines r\'esolutions, qui conduisent
\`a des produits de
convolution de la forme $\calA_{\mu_1}*\cdots*\calA_{\mu_n}$ o\`u chaque
$\mu_i$ est minuscule ou
quasi-minuscule.
L'id\'ee essentielle de la d\'emonstration est que chaque $\calA_\lambda$
appara\^{\i}t
comme facteur direct (avec une certaine multiplicit\'e) d'un tel produit.
Cette assertion
est ramen\'ee \`a un \'enonc\'e combinatoire, dont nous donnons deux
preuves diff\'erentes dans la section 10.
L'une repose sur un lemme simple sur les syst\`emes des racines,
l'autre est bas\'ee sur la th\'eorie des repr\'esentations et le mod\`ele
des chemins de Littelmann.

Arm\'es de la connaissance explicite des cas minuscule et quasi-mi\-nus\-cule,
et des r\'esultats de la section 9,
on peut alors d\'emontrer les th\'eor\`emes 3.1 et 3.2 en suivant
l'argument de \cite{Ngo}. Ceci est le contenu de la section 11.

Nos r\'esultats ont \'et\'e expos\'es au s\'eminaire Formes automorphes de
l'Uni\-ver\-sit\'e Paris 7 en F\'evrier 1999,
et au Number theory seminar
du Max Planck Institut fuer Mathematik en Juin 1999. En Juillet 1999,
Frenkel, Gaitsgory et Vilonen ont annonc\'e une autre
d\'emonstration de la conjecture de \cite{FGKV}, par une voie diff\'erente
et ind\'ependante \cite{FGV}.

Pendant la r\'edaction de ce travail, B.C. Ng\^o
a b\'en\'efici\'e de l'hospita\-li\-t\'e du
Max Planck Institut fuer Mathematik.
Il remercie G. Laumon et M. Rapoport pour
d'utiles discussions sur le sujet de cet article.

Nous remercions le rapporteur pour sa lecture attentive du manuscrit.

\section{Notations}

Soit $k$ un corps fini \`a $q$ \'el\'ements et de caract\'eristique $p$.
Notons ${\bar k}$ sa cl\^oture alg\'ebrique.
Nous supposons, pour des raisons de commodit\'e, que $G$ est un groupe
alg\'ebrique semi-simple d\'eploy\'e sur $k$, la
g\'en\'eralisation aux groupes r\'eductifs \'etant \'evidente. Soient $T$ un
tore maximal d\'eploy\'e de $G$, et $B,B^-$ deux
sous-groupes de Borel tels que $B\cap B^- = T$. On note $U$ (resp. $U^-$)
le radical unipotent de $B$ (resp. $B^-$). On pose
${\mathfrak g} = \Lie(G)$ et ${\mathfrak h} = \Lie(T)$.

On note $\langle\phantom{a},\phantom{a}\rangle$ l'accouplement naturel entre
$X := \Hom(T,\Gm)$ et $X^\smv := \Hom(\Gm,T)$.
Soient $R \subset X$ le syst\`eme de racines associ\'e \`a $(G,T)$, $R_+$
(resp $R_-$) l'ensemble des racines correspondant \`a $B$ (resp. $B^-$) et
$\Delta=\{\alpha_1,\ldots,\alpha_r\}$
l'ensemble des racines simples. Pour tout $\alpha\in R$, on
d\'esigne par $U_\alpha$ le sous-groupe radiciel de $G$
correspondant \`a la racine $\alpha$. Soit $R^\smv\subset X^\smv$ le
syst\`eme de racines dual muni de la bijection
$R\to R^\smv$, $\alpha\mapsto\alpha^\smv$. Notons $R^\smv_+$ l'ensemble des
coracines positives.
Soit $W$ le groupe de Weyl de $(G,T)$.

Soit $\rho=(1/2)\sum_{\alpha\in R_+}\alpha$ la demi-somme des racines
positives.
Pour toute racine simple $\alpha\in\Delta$, on a
$\langle\rho,\alpha^\smv\rangle=1$.

On note $Q^\smv$ (resp. $Q^\smv_+$) le sous-groupe (resp. le
sous-mono\"{\i}de) de $X^\smv$ engendr\'e par $R^\smv$ (resp. $R^\smv_+$).
On d\'esigne par $X^\smv_+$ le c\^one des cocaract\`eres dominants :
$$
X_+^\smv=\{\lambda\in X^\smv\mid \langle\alpha,\lambda\rangle\geq 0,\;
\forall \alpha\in R_+\}.
$$

On consid\`ere l'ordre partiel sur $X^\smv$ d\'efini comme suit :
$\nu\geq\nu'$ si et seulement si $\nu-\nu'\in Q^\smv_+$.

On note $G^\smv$ le groupe dual, consid\'er\'e sur $\Ql$; il est muni des
sous-groupes
$T^\smv \subset B^\smv$. Pour tout $\lambda\in X^\smv_+$, on note
$$
\Omega(\lambda) = \{\nu\in X^\smv\mid \forall w\in W, w\nu\leq\lambda\};
$$
c'est l'ensemble des poids de $T^\smv$ dans
$V(\lambda)$, le $G^\smv$-module simple, sur $\Ql$, de plus haut poids
$\lambda$
(voir, par exemple, \cite[Chap.VIII, Ex.7.1]{Bourbaki} ou
\cite[Prop. 23.1]{Humphreys}).

\medskip
On notera $M$ l'ensemble des \'el\'ements minimaux de $X^\smv_+\setminus
\{0\}$.

\begin{lemme}
Soit $\mu\in M$. On a l'une des alternatives suivantes.
\begin{enumerate}
\item
Si $\langle \alpha,\mu \rangle \in \{0, \pm 1\}$ pour tout
$\alpha\in R$, alors $\mu$ est un \'el\'ement minimal de $X^\smv_+$ et l'on
a\,
$\Omega(\mu) = W\mu$.
Dans ce cas, on dira que $\mu$ est un cocaract\`ere {\it minuscule}.
\item
Sinon, il existe une unique racine $\gamma$ telle que $\langle\gamma, \mu
\rangle \geq 2$; c'est
une racine positive maximale, et l'on a $\mu = \gamma^\smv$ et\,
$\Omega(\mu) = W\mu \cup \{0\}$.
Dans ce cas, on dira que $\mu$ est {\it quasi-minuscule}.
\end{enumerate}
\end{lemme}

\dem Compte-tenu des r\'ef\'erences cit\'ees avant le lemme, la premi\`ere
assertion
r\'esulte de \cite[Chap.VI, Ex.1.24]{Bourbaki}. Voyons la seconde.

Soit $\gamma\in R$ tel que $\langle\gamma, \mu \rangle \geq 2$.
D'apr\`es \cite[Chap.VI, Ex.1.23]{Bourbaki} ou \cite[Prop. 23.1]{Humphreys},
$\mu-\gamma^\smv$ est $W$-conjugu\'e \`a un poids
dominant $\leq \mu$, donc \`a $0$ ou $\mu$ (puisque $\mu\in M$). Or, comme
$\langle\gamma, \mu \rangle \geq 2$,
on voit facilement que la norme de $\mu-\gamma^\smv$ (relativement \`a un
produit scalaire $W$-in\-va\-riant)
est strictement inf\'erieure \`a celle de $\mu$. On en d\'eduit que $\mu =
\gamma^\vee$, et que
$\gamma$ est l'unique racine telle que $\langle\gamma, \mu \rangle \geq 2$.
Soient $R_\gamma$
(resp. $R^\smv_{\gamma^\smv}$) le sous-syst\`eme de racines irr\'eductible
de $R$
(resp. $R^\smv$) contenant $\gamma$ (resp. $\gamma^\smv$); on dira que les
\'el\'ements de
$R^\smv_{\gamma^\smv}$
de longueur minimale sont des coracines courtes.
Il est alors bien connu que l'\'egalit\'e $\Omega(\gamma^\smv) =
W\gamma^\smv \cup \{0\}$
entra\^{\i}ne que $\gamma$ est l'unique racine maximale de
$R_\gamma$; de fa\c con
\'equivalente,
$\gamma^\smv$ est l'unique coracine courte dominante  de
$R^\smv_{\gamma^\smv}$ (\cf \cite[Chap.VIII, 7.22]{Bourbaki}.
\findem

\smallskip{\sc Remarque.} L'usage du mot minuscule est ici plus g\'en\'eral
que celui de
\cite{Bourbaki}, qui se limite au cas o\`u $R$ est irr\'eductible, auquel
cas un copoids minuscule
est n\'ecessairement un copoids fondamental.

\medskip
Soient $\calO=k[[\varpi]]$ l'anneau des s\'eries formelles en une variable
$\varpi$ et $F=k(\!(\varpi)\!)$ son corps des fractions.

Pour chaque \'el\'ement $\nu\in X^\smv$, on note $\varpi^\nu\in T(F)$
l'image par le coca\-ract\`ere $\nu$,
de l'uniformisante $\varpi\in F^\times$. On rappelle (voir
\cite[\S3.5]{Cartier} et \cite{MacDonald}) les d\'ecompositions de
Cartan et d'Iwasawa
$$
\begin{array}{rcl} G(F) &=&\coprod_{\lambda\in
X^\smv_+}G(\calO)\varpi^\lambda G(\calO) ; \\
G(F) &=& \coprod_{\nu\in X^\smv} U(F)\varpi^\nu G(\calO) .
\end{array}
$$

\noindent{\it Convention.\ } Sauf mention expresse du contraire,
quand on parle de points de $k$-sch\'emas, il s'agira des points
\`a valeurs dans une $k$-alg\`ebre arbitraire. Par {\it stratification}
d'un $k$-sch\'ema $X$, on entend la donn\'ee
des sous-sch\'emas localement ferm\'es $X_\alpha$ de $X$,
deux \`a deux disjoints, tels que $X=\cup X_\alpha$.

\section{La Grassmannienne affine}

Rappelons la construction de la Grassmannienne affine $\Gr$, tir\'ee de
\cite{Beauville-Laszlo} et \cite{Laszlo-Sorger}. Comme
dans \loc, appelons $k$-espaces, resp. $k$-groupes, les faisceaux
d'en\-sem\-bles, resp. de groupes, sur la cat\'egorie des
$k$-alg\`ebres munie de la topologie fid\`element plate et de
pr\'esentation finie.

Consid\'erons le $k$-groupe $LG$ et
son $k$-sous-groupe $L^{\geq 0}G$,
qui associent \`a chaque $k$-alg\`ebre $R$, le groupe $G(R(\!(\varpi)\!))$
et le sous-groupe $G(R[[\vp]])$.
Ces constructions s'appliquent aussi aux sous-groupes $T$ et $U$ de $G$,
\`a la place de $G$.

Il est clair que $L^{\geq 0}G$ est repr\'esent\'e par le sch\'ema en groupes
limite projective des sch\'emas en groupes de type fini $R\mapsto
G(R[[\vp]]/(\vp^n))$.
Pour d\'efinir une structure d'ind-sch\'ema
sur $LG$, choisissons une repr\'e\-sen\-ta\-tion fid\`ele $\rho:G\rightarrow SL(V)$.
Notons $L^{(N)}G(R)$ l'ensemble des $g\in LG(R)$ tel
que l'ordre des p\^oles de $\rho(g)$ et de $\rho(g^{-1})$ n'exc\`ede pas
$N$. D'apr\`es \loc, $L^{(N)}G$ est repr\'esentable
par un sch\'ema, et le faisceau $\Gr$ associ\'e au pr\'efaisceau
$R\mapsto G(R(\!(\varpi)\!))/G(R[[\vp]])$ est une limite inductive de
sch\'emas projectifs $\Gr^{(N)}=L^{(N)}G/L^{\geq 0}G$.

Notons $L^{\leq 0} G$ le $k$-groupe $R\mapsto G(R[\vp^{-1}])$ et $L^{<0}G$
le noyau du morphisme $L^{\leq 0} G\rightarrow G$
d\'efini par $\vp^{-1}\mapsto 0$. Ce sont des $k$-sous-groupes de $LG$.

D'apr\`es \cite[Prop. 1.11]{Beauville-Laszlo} et \cite[Prop.
4.6]{Laszlo-Sorger},
on a alors le lemme suivant. Dans \loc, $G$ est
suppos\'e simplement connexe et $k={\mathbb C}$, mais la d\'emons\-tration
 s'\'etend au cas g\'en\'eral. Ce r\'esultat d\'ecoule aussi d'un
th\'eor\`eme de Ramanathan \cite{Ramanathan}.

\begin{lemme}
Le morphisme de multiplication
$$L^{<0}G\times L^{\geq 0}G\rightarrow LG$$ est une immersion ouverte.
\end{lemme}

Identifions $L^{<0}G$ \`a l'ouvert $L^{<0}G e_0$, o\`u $e_0$ d\'esigne le
point base de $\Gr$. La Grassmannienne affine $\Gr$
est recouverte par les ouverts translat\'es $gL^{<0}Ge_0$ au-dessus desquels
la fibration $LG\rightarrow\Gr$ est triviale.
Ces ouverts trivialisants sont utiles pour \'etudier
de mani\`ere plus explicite la g\'eom\'etrie locale de
$\Gr$. Par exemple, $L^{<0}G$ n'est pas r\'eduit en g\'en\'eral si bien que
$\Gr$ ne l'est pas non plus.

Le groupe $L^{\geq 0}G$ agit naturellement sur $\Gr$.
Pour tout $\lambda\in X^\smv$,
notons $e_\lambda$ le point $\vp^\lambda e_0$ de $\Gr$. Pour $\lambda\in
X^\smv_+$, notons $\Gr_\lambda$ l'orbite de
$L^{\geq 0}G$ passant par $e_\lambda$. Notons ${\bar\Gr}_\lambda$
l'adh\'erence de $\Gr_\lambda$. Introduisons aussi
les sous-groupes
$L^{\geq\lambda}G :=\vp^{\lambda}L^{\geq 0}G\vp^{-\lambda}$ et
$L^{<\lambda}G :=\vp^{\lambda}L^{< 0}G\vp^{-\lambda}$.

Notons $J$ l'image inverse du radical unipotent $U$ de $B$ par
l'homo\-mor\-phisme $L^{\geq 0}G\rta G$ d\'efini par $\varpi\mapsto 0$ ; c'est
une limite projective de groupes unipotents. Posons
$J^{\geq\lambda}=J\cap L^{\geq\lambda}G$ et $J^\lambda=J\cap
L^{<\lambda}G$.

Quelques soient $\alpha\in R$ et $i\in \ZZ$, on d\'esigne par
$U_{\alpha,i}$ l'image de l'homomor\-phisme
${\mathbb G}_a\rta LG$ d\'efini par $x\mapsto U_\alpha(\vp^i x)$.
La multiplication d\'efinit un isomorphisme
$$
\prod_{\alpha\in R_+,\langle\alpha,\lambda\rangle>0}
\prod_{i=0}^{\langle\alpha,\lambda\rangle-1} U_{\alpha,i}\rta J^\lambda
$$
(en choisissant un ordre total sur l'ensemble des facteurs). En
particulier, $J^\lambda$ est isomorphe \`a l'espace affine
de dimen\-sion $2\langle\rho,\lambda\rangle$.

\begin{lemme}
Le morphisme naturel $J^\lambda\rta \Gr_\lambda$ d\'efini par $j\mapsto
je_\lambda$ est une immersion ouverte. \end{lemme}

\dem Il est clair que la multiplication induit un isomorphisme
$J^\lambda\times J^{\geq\lambda}\rta J$.
Il est aussi clair que la multiplication induit une immersion ouverte
$J\times B^-\rta L^{\geq 0} G$. Par ailleurs,
$J^{\geq\lambda}$ et $B^-$ sont des sous-groupes de $L^{\geq\lambda}G$ qui
fixent $e_\lambda$. Le lemme s'en d\'eduit.
\findem

Il r\'esulte du lemme que l'orbite $\Gr_\lambda$ est lisse, irr\'eductible et
de dimension $2\langle\rho,\lambda\rangle$. Elle est incluse dans un
$\Gr^{(N)}$ pour $N\in{\mathbb N}$
assez grand si bien que son adh\'erence
${\bar\Gr}_\lambda$ est un sch\'ema projectif, irr\'eductible et stable par
l'action de $L^{\geq 0}G$.
Il est bien connu, voir \cite[\S11]{Lusztig}, que ${\bar\Gr}_\lambda$
est la r\'eunion des orbites ${\Gr}_{\lambda'}$ avec $\lambda'\leq\lambda$.
En particulier, si $\mu\in X^\smv_+$ est nul ou bien minuscule, l'orbite
$\Gr_\mu$
est un sch\'ema projectif
lisse.

Notons $L^{>0}G$ le noyau de l'homomorphisme $L^{\geq 0}G\rta G$;
c'est une limite projective de groupes unipotents. Il est clair
que pour tout $\lambda\in X^\smv_+$, le morphisme
$$
(L^{>0}G\cap L^{\geq\lambda}G)\times (L^{>0}G\cap L^{<\lambda}G)\rta L^{>0}G$$
est un isomorphisme et que
$$L^{>0}G\cap L^{<\lambda}G
=\prod_{\alpha\in R_+,\langle\alpha,\lambda\rangle>1}
\prod_{i=1}^{\langle\alpha,\lambda\rangle-1} U_{\alpha,i}.$$

Soit $P_\lambda$ le sous-groupe parabolique de $G$ engendr\'e par $B^-$ et
par les sous-groupes radiciels $U_\alpha$ avec
$\langle\alpha,\lambda\rangle=0$. Le groupe de Weyl de $P_\lambda$ est
\'egal au stabilisateur $W_\lambda$
de $\lambda$. Notons $U^+_\lambda$ le radical unipotent du parabolique
oppos\'e \`a $P_\lambda$. Il est clair que
$P_\lambda \subset L^{\geq \lambda}G$ et que
$J^\lambda= U^+_\lambda\ltimes(L^{>0}G\cap L^{<\lambda}G)$.

\begin{lemme}
On a
$$L^{\geq 0}G\cap L^{\geq\lambda}G=P_\lambda\ltimes (L^{>0}G\cap
L^{\geq\lambda}G).$$
En particulier, le groupe
$L^{\geq 0}G\cap L^{\geq\lambda}G$ est g\'eom\'etrique\-ment connexe,
et l'on a $G\cap L^{\geq\lambda}G = P_\lambda$.  \end{lemme}

\dem
Il suffit de d\'emontrer que le morphisme de multiplication
$$(L^{>0}G\cap L^{\geq\lambda}G)\times P_\lambda \rta L^{\geq 0}G\cap
L^{\geq\lambda}G$$
est un isomorphisme.

Soit $g$ un point de $L^{\geq 0}G$ qui s'\'ecrit sous la forme
$g=g^+ g^- u w p$ o\`u $g^+\in (L^{>0}G\cap L^{\geq\lambda}G)$,
o\`u $g^-\in (L^{>0}G\cap L^{<\lambda}G)$, o\`u $u\in (U \cap wU^+_\lambda
w^{-1})$,
o\`u $p\in P_\lambda$
et o\`u $w$ est de longueur minimale dans sa classe $wW_\lambda$.

Supposons que $g\in L^{\geq\lambda}G$.
Puisque $g^+$ et $p$ appartiennent d\'ej\`a \`a ce groupe, il en est de
m\^eme de
$g^- uw$. On a donc
$$\vp^{-\lambda} g^- uw\vp^{\lambda}=(\vp^{-\lambda} g^- u\vp^{\lambda})
\vp^{w\lambda-\lambda}w
\in L^{\geq 0}G.$$
Puisque $w$ appartient \`a $L^{\geq 0}G$,
$g^- u$ appartient \`a $LU$
et $\vp^{w\lambda-\lambda}\in LT$ et compte tenu de la d\'ecompostion
d'Iwahori \cite[3.5]{Cartier}, $\vp^{-\lambda} g^-
u\vp^{\lambda}$ et $\vp^{w\lambda-\lambda}$ appartiennent,
tous les deux, \`a $L^{\geq 0}G$. Par ailleurs, $\vp^{w\lambda-\lambda}$
appartient \`a $L^{\geq 0}G$ si et seulement si $w\lambda=\lambda$.
Comme $w$ est de
longueur minimale dans sa classe $wW_\lambda$, il vient $w=1$. Donc, on a
$u\in U^+_\lambda$. Compte tenu de la
d\'ecomposition $J^\lambda= U^+_\lambda\ltimes(L^{>0}G\cap L^{<\lambda}G)$
et du fait que $J^\lambda\cap L^{\geq 0}G=\{1\}$, on obtient
$g^-=1$ et $u=1$.
\findem

Soit $\ell$ un nombre premier diff\'erent de la caract\'eristique de $k$.
Pour tout $\lambda\in X^\smv_+$, notons $\calA_\lambda$ le complexe
d'intersection
$\ell$-adique de ${\bar\Gr}_\lambda$. D'apr\`es le lemme pr\'ec\'edent,
le stabilisateur de chaque $e_\lambda$ dans $L^{\geq 0}G$ est
g\'eom\'etrique\-ment connexe.   
Par cons\'equent, tout faisceau pervers sur $\Gr$, g\'eom\'etriquement
irr\'eductible, $L^{\geq 0}G$-\'equivariant et dont le support est un
sch\'ema de type fini, est isomorphe \`a ${\calA}_\lambda$ pour
un certain $\lambda\in X^\smv_+$.

A la suite de
Lusztig, Ginzburg, Mirkovic et Vilonen, voir \cite{Lusztig,Lusztig97},
\cite{Ginzburg} et \cite{Mirkovic-Vilonen},
on va d\'efinir
le produit de convolution $\calA_{\lambda_1}*\calA_{\lambda_2}$ pour
$\lambda_1,\lambda_2\in X^\smv_+$ comme suit. Consid\'erons
les morphismes
$$\Gr\times\Gr\ \hflg{\pi_1}{}\ LG\times\Gr \
\hfld{\pi_2}{}\ \Gr\times\Gr$$
d\'efinis par $\pi_1(g,x)=(ge_0,x)$ et
$\pi_2(g,x)=(ge_0,gx)$. Le morphisme $\pi_1$ est le morphisme quotient pour
l'action de $L^{\geq 0}G$ sur $LG\times\Gr$
d\'efinie par
$$\alpha_1(h)(g,x)=(gh^{-1},x).$$
Le morphisme $\pi_2$ est le
morphisme quotient pour l'action de
$L^{\geq 0}G$ sur $LG\times\Gr$ d\'efinie par
$$\alpha_2(h)(g,x)=(gh^{-1},hx).$$ Pour tous $\lambda_1,\lambda_2\in
X^\smv_+$,
notons ${\bar\Gr}_{\lambda_1}\tildetimes{\bar\Gr}_{\lambda_2}$ le
quotient de
$\pi_1^{-1}({\bar\Gr}_{\lambda_1}\times{\bar\Gr}_{\lambda_2})$ par
$\alpha_2(L^{\geq 0}G)$. L'existence de ce quotient est assur\'ee par la
locale trivialit\'e du morphisme $LG\rta\Gr$.
Plus pr\'ecis\'ement, au-dessus des ouverts de ${\bar\Gr}_{\lambda_1}$ de
la forme $g L^{<0}G e_0\cap{\bar\Gr}_{\lambda_1}$, les sch\'emas
${\bar\Gr}_{\lambda_1}\,{\tilde\times}\,{\bar\Gr}_{\lambda_2}$ et
${\bar\Gr}_{\lambda_1}\times{\bar\Gr}_{\lambda_2}$ sont isomorphes.
De plus, ces isomorphismes sont clairement compatibles avec
la stratification de
${\bar\Gr}_{\lambda_1}\times{\bar\Gr}_{\lambda_2}$ par les sous-sch\'emas
localement ferm\'es
$\Gr_{\lambda'_1}\times\Gr_{\lambda'_2}$ et celle de
${\bar\Gr}_{\lambda_1}\,{\tilde\times}\,{\bar\Gr}_{\lambda_2}$ par les
sous-sch\'emas localement ferm\'es
$${\Gr}_{\lambda'_1}\,{\tilde\times}\,\Gr_{\lambda'_2}
=\pi_1^{-1}({\Gr}_{\lambda'_1}\times{\Gr}_{\lambda'_2})/\alpha_2
(L^{\geq 0}G),$$
avec $\lambda_1'\leq\lambda_1$ et $\lambda'_2\leq\lambda_2$. La projection
sur le second facteur d\'efinit un morphisme
$$m:{\bar\Gr}_{\lambda_1}\,{\tilde\times}\,{\bar\Gr}_{\lambda_2} \rta
{\bar\Gr}_{\lambda_1+\lambda_2}.$$
On pose
$$\calA_{\lambda_1}*\calA_{\lambda_2}={\rm R}m_*
(\calA_{\lambda_1}\,{\tilde\boxtimes}\,\calA_{\lambda_2}),$$ o\`u
$\calA_{\lambda_1}\,{\tilde\boxtimes}\,\calA_{\lambda_2}$ d\'esigne le
complexe d'intersection de 
${\bar\Gr}_{\lambda_1}\tildetimes{\bar\Gr}_{\lambda_2}$.

La construction pr\'ec\'edente, g\'en\'eralis\'ee de la mani\`ere \'evidente,
permet de d\'efinir le produit de convolution it\'er\'e
$$
\calA_{\lambda_1}*\cdots*\calA_{\lambda_n}
$$
pour tous
$\lambda_1,\ldots,\lambda_n\in X^\smv_+$. D'apr\`es \cite{Ginzburg} et
\cite{Mirkovic-Vilonen}, ce produit de convolution est encore un faisceau
pervers, somme directe, avec multiplicit\'es, de faisceaux
pervers $\calA_\lambda$ avec $\lambda\leq\lambda_1+\cdots+\lambda_n$. Nous
n'utiliserons ce r\'esultat que dans le cas o\`u les $\lambda_i$
appartiennent \`a l'ensemble $M$. Nous proposons une d\'emonstration
simple dans ce cas et montrons comment le cas g\'en\'eral peut, en fait, se
d\'eduire de celui-ci.

\section{Les \'enonc\'es principaux}
Rappelons que $U$ d\'esigne le radical unipotent du sous-groupe de Borel
$B$ associ\'e \`a $R_+$. On d\'efinit de fa\c con analogue pour $U$, \`a la
place de $G$,
les groupes de lacets $LU$, $L^{\geq 0}U=LU\cap L^{\geq 0}G$ et
$L^{<0}U=LU\cap L^{<0}G$. Pour tout $\nu\in X^\smv$, on note aussi $L^{\geq
\nu}U=\vp^\nu L^{\geq 0}U\vp^{-\nu}$ et
$L^{<\nu}U=\vp^{\nu}L^{<0}U\vp^{-\nu}$. La multiplication induit un
isomorphisme
$$L^{\geq\nu}U\times L^{<\nu}U\rta LU.$$

Pour tout $\nu\in X^\smv$, $L^{<\nu}U$ est un sous-groupe ferm\'e de
$L^{<\nu}G$
si bien qu'on peut identifier $L^{<\nu}U e_0$ \`a un ferm\'e, not\'e
$S_\nu$, de l'ouvert $\vp^\nu L^{<0}Ge_0$ de $\Gr$. En
particulier, pour tout $\lambda\in X^\smv_+$ et tout $\nu\in X^\smv$,
$S_\nu\cap{\bar\Gr}_\lambda$ est un sous-sch\'ema
localement ferm\'e, \'eventuellement vide, de ${\bar\Gr}_\lambda$.
D'apr\`es la d\'ecomposition d'Iwasawa, ces intersections
$S_\nu\cap{\bar\Gr}_\lambda$ forment une stratification de ${\bar\Gr}_\lambda$.

Nous donnerons une nouvelle d\'emonstration du th\'eor\`eme suivant, d\^u
\`a Mirkovic et Vilonen dans le cas $k={\mathbb C}$
(\cite{Mirkovic-Vilonen}).

\begin{thr} Quelques soient $\lambda\in X^\smv_+$ et $\nu\in X^\smv$,  
le complexe ${\rm R}\Gamma_c(S_\nu,\calA_\lambda)$ est concentr\'e en
degr\'e $2\langle\rho,\nu\rangle$.
De plus, l'endo\-morphisme ${\rm Fr}_q$ agit dans ${\rm
H}_c^{2\langle\rho,\nu\rangle}(S_\nu,\calA_\lambda)$ comme
$q^{\langle\rho,\nu\rangle}$.
\end{thr}

Dans l'\'enonc\'e pr\'ec\'edent on a \'ecrit 
${\rm R}\Gamma_c(S_\nu,\calA_\lambda)$
\`a la place de
$${\rm R}\Gamma_c((S_\nu\cap{\bar\Gr}_\lambda)\otimes_k{\bar k},
\calA_\lambda),$$
pour all\'eger la notation. 
On utilisera syst\'ematiquement cette notation all\'eg\'ee dans la suite,
cet abus de notation
ne causant aucune ambiguit\'e.

Soient $\nu\in X^\smv_+$ et $\nu'\in -X^\smv_+$. En choisissant un ordre
total sur les racines positives, on a un isomorphisme
$$\prod_{\alpha\in R^+}
\prod_{\langle\alpha,\nu'\rangle\leq i<\langle\alpha,\nu\rangle}
U_{\alpha,i}=L^{<\nu}U\cap L^{\geq\nu'}U.$$
Pour $\nu$ fix\'e et pour $\nu'$ de plus en plus anti-dominant, ces groupes
forment un syst\`eme inductif dont la limite est $L^{<\nu}U$.

Pour toute racine simple $\alpha\in\Delta$, notons $u_{\alpha,i}$ la
projection sur le facteur $U_{\alpha,i}$ et $$h:L^{<\nu}U\cap
L^{\geq\nu'}U\rta {\mathbb G}_a$$ le morphisme
$h(x)=\sum_{\alpha\in\Delta}u_{\alpha,-1}(x)$. Ce morphisme est visiblement
compa\-tible aux fl\`eches de transition et induit sur la limite inductive un
morphisme $h:L^{<\nu}U\rta{\mathbb G}_a$, pour tout
$\nu$ dominant. Compte tenu de l'isomorphisme $L^{<\nu}U \times
L^{\geq\nu}U\rta LU$, $u^-u^+\mapsto u$, on peut d\'efinir un
morphisme, not\'e aussi $h: LU\rta{\mathbb G}_a$, par la relation
$h(u^-u^+)=h(u^-)$. Ce morphisme ne d\'epend pas du $\nu$
dominant choisi.
On en d\'eduit un morphisme, not\'e encore $h:S_\nu\rta{\mathbb G}_a$.

Fixons un caract\`ere additif
non trivial $\psi:k\rta\Ql^\times$ et notons ${\mathcal L}_\psi$ le
faisceau d'Artin-Schreier sur ${\mathbb G}_a$ associ\'e \`a $\psi$. Le
caract\`ere $\theta:U(F)\rta\Ql$ consid\'er\'e dans l'introduction est le
caract\`ere $x\mapsto \psi(h(x))$.

L'\'enonc\'e suivant a \'et\'e conjectur\'e par Frenkel, Gaitsgory, Kazhdan
et Vilonen \cite{FGKV}.

\begin{thr}
Pour $\nu\not=\lambda$ dans $X^\smv_+$, le complexe ${\rm
R}\Gamma_c(S_\nu,\calA_\lambda\otimes h^*{\mathcal L}_\psi)$ est nul.
Pour $\nu=\lambda$, ce complexe est isomorphe \`a $\Ql$, muni de l'action
de Frobenius agissant par $q^{\langle\rho,\lambda\rangle}$, plac\'e en
degr\'e $2\langle\rho,\lambda\rangle$. \end{thr}

Ces r\'esultats impliquent les \'enonc\'es sur les termes constants et les
coefficients de Fourier
mentionn\'es dans l'introduction via le dictionnaire faisceaux-fonctions de
Grothendieck \cite{Grothendieck}.

Nous exposerons les d\'emonstrations de ces deux th\'eor\`emes en
parall\`ele dans la suite de l'article.

\section{L'action du tore $T$}

Le tore $T$ normalise les sous-groupes $L^{\geq 0}G$, $L^{<0}G$,
$L^{<\nu}U$ ... de $LG$ si bien qu'il agit sur tous les objets
g\'eom\'etriques
qu'on a consid\'er\'es dans les deux sections pr\'ec\'edentes. Cette action
fournit un outil pr\'ecieux pour \'etudier leur g\'eom\'etrie.
Choisissons une fois pour toutes un cocaract\`ere strictement dominant
$\phi:{\mathbb G}_m\rta T$. On entendra par l'action $\phi({\mathbb G}_m)$,
l'action restreinte de $T$ \`a ${\mathbb G}_m$ via ce cocaract\`ere.

\begin{lemme}
Pour tout $\nu\in X^\smv$, le point $e_\nu$ est le seul point fixe de
l'action $\phi({\mathbb G}_m)$ sur $S_\nu$. De plus, c'est un point fixe
attractif.
\end{lemme}

\dem
Tout $x\in L^{<\nu}U({\bar k})$ est de la forme
$$x=\prod_{\alpha\in R_+}\prod_{i<\langle\alpha,\nu\rangle} U_{\alpha,i}
(x_{\alpha,i}),
$$
o\`u les $x_{\alpha,i}\in{\bar k}$ sont nuls sauf pour un nombre fini
d'entre eux.
Alors, pour tout $z\in{\bar k}^\times$, on a
$$
\phi(z) xe_\nu=\prod_{\alpha\in R^+}\prod_{i<\langle\alpha,\nu\rangle}
U_{\alpha,i}
(z^{\langle\alpha,\phi\rangle}x_{\alpha,i})e_\nu.$$ Le lemme r\'esulte donc
de l'hypoth\`ese que pour tout $\alpha\in R_+$,
l'entier $\langle\alpha,\phi\rangle$ est strictement positif. \findem

Ce lemme montre que les $e_\nu$ sont les seuls points fixes de l'action
$\phi({\mathbb G}_m)$ sur $\Gr$.
De plus, il entra\^{\i}ne l'\'enonc\'e suivant, qui pr\'ecise
\cite[2.6.11(3)]{MacDonald} et est certainement bien connu.

\begin{lemme}
Si l'intersection $S_\nu\cap{\bar\Gr}_\lambda$ n'est pas vide, $\nu$
appartient \`a $\Omega(\lambda)$.
\end{lemme}

\dem
Si un point $xe_\nu$, avec $x\in L^{<\nu}U({\bar k})$, appartient \`a
${\bar\Gr}_\lambda({\bar k})$,
toute l'orbite de $\phi({\mathbb G}_m)$ passant par ce point y appartient
aussi. Puisque le point fixe $e_\nu$ appartient \`a
l'adh\'erence de cette orbite et que ${\bar\Gr}_\lambda$ est propre,
$e_\nu$ appartient \`a ${\bar\Gr}_\lambda$ d'o\`u
$\nu\in\Omega(\lambda)$. \findem

Signalons au passage l'\'enonc\'e suivant qui ne servira pas dans la suite
de l'article. Cet \'enonc\'e a \'et\'e d\'ecouvert lors
d'une conversation que l'un de nous a eu avec M. Rapoport.

\begin{proposition}
La caract\'eristique d'Euler-Poincar\'e
$\chi_c(S_\nu\cap\Gr_\lambda)$ est \'egale \`a $1$ si $\nu$ est conjugu\'e \`a
$\lambda$ par un \'el\'ement de $W$ et \`a $0$ sinon.
\end{proposition}

\dem
Dans le premier cas, $S_\nu\cap\Gr_\lambda$ contient un unique point fixe
$e_\nu$ de l'action $\phi({\mathbb G}_m)$. Dans le second cas, le groupe
multiplicatif $\phi({\mathbb G}_m)$ agit sans point fixe. La proposition
r\'esulte donc d'un th\'eor\`eme de Bialynicki-Birula \cite[cor.
2]{Bialynicki-Birula}.
\findem

Cet \'enonc\'e peut \^etre consid\'er\'e comme une interpr\'etation
g\'eom\'etri\-que d'un r\'esultat de Lusztig \cite[6.1]{Lusztig}.
Reprenons les notations de l'intro\-duc\-tion.
Soit $C_\lambda$ l'\'el\'ement de l'alg\`ebre
de Hecke $\calH$ d\'efini par
$$C_\lambda=(-1)^{2\langle\rho,\lambda\rangle}
q^{-\langle\rho,\lambda\rangle}{\mathbb I}_\lambda,$$
o\`u ${\mathbb I}_\lambda$ est la fonction caract\'eristique
de $K\varpi^\lambda K$. On sait que
$$(C_\lambda)=(K_{\lambda,\mu}(q))^{-1}(A_\lambda), $$
o\`u $(K_{\lambda,\mu}(q))$ est la matrice triangulaire
form\'ee des polyn\^omes de Kazhdan et Lusztig. Les termes
constants normalis\'es
$$(-1)^{2\langle\rho,\nu\rangle}q^{-\langle\rho,\nu\rangle}
\int_{U(F)}C_\lambda(x\vp^\mu) dx$$
sont donc calcul\'es par la matrice
$$(K_{\lambda,\mu}(q))^{-1}(m_\lambda(\mu))$$
d'apr\`es le th\'eor\`eme de Lusztig-Kato.
Compte tenu de la proposition pr\'ec\'e\-den\-te
(et aussi de 3.1), on obtient, en sp\'ecialisant $q\mapsto1$,
$$(K_{\lambda,\mu}(1))^{-1}(m_\lambda(\mu))={\rm Id}, $$
d'o\`u $K_{\lambda,\mu}(1)=m_\lambda(\mu)$.

\section{Les intersections $S_{w\lambda}\cap{\bar\Gr}_\lambda$}

Pour $\lambda\in X^\smv_+$, on a consid\'er\'e dans la section 2 le groupe

$$J^\lambda=\prod_{\alpha\in R_+}
\prod_{i=0}^{\langle\alpha,\lambda\rangle-1}U_{\alpha,i}$$ qui est
manifestement un sous-groupe de $L^{\geq 0}U$.
On a aussi d\'emontr\'e que le morphisme
$J^{\lambda}\rta{\bar\Gr}_{\lambda}$ d\'efini par $j\mapsto je_\lambda$ est
une immersion
ouverte.

\begin{lemme}Soit $\lambda\in X^\smv_+$. Le morphisme $j\mapsto je_\lambda$
induit un isomorphisme de
$J^\lambda$ sur l'ouvert $\vp^\lambda L^{<0}Ge_0\cap{\bar\Gr}_\lambda$
de ${\bar\Gr}_\lambda$.
\end{lemme}

\dem
L'image de $J^\lambda$ est contenue dans
$\vp^\lambda L^{<0}Ge_0\cap{\bar\Gr}_\lambda$.
D'apr\`es le lemme 2.2, elle est en fait un ouvert dense de $\vp^\lambda
L^{<0}Ge_0\cap{\bar\Gr}_\lambda$.
Par cons\'equent, $\vp^{-\lambda}J^\lambda\vp^{\lambda}$ est un ouvert
dense dans l'image inverse de
$\vp^\lambda L^{<0}Ge_0\cap{\bar\Gr}_\lambda$
par l'isomorphisme
$$L^{<0}G\rta \vp^\lambda L^{<0}Ge_0.$$
Or, $\vp^{-\lambda}J^\lambda\vp^{\lambda}$ est un sous-groupe ferm\'e de
$L^{<0}G$ et le lemme s'en d\'eduit.
\findem

\begin{lemme} Soit $\lambda\in X^\smv_+$.
Pour tout $w\in W$, le morphisme
$$wJ^\lambda w^{-1}\cap LU\rta S_{w\lambda}\cap{\bar\Gr}_\lambda$$ d\'efini
par $j\rta j e_{w\lambda}$ est un isomorphisme.
Par cons\'equent, $S_{w\lambda}\cap{\bar\Gr}_\lambda$
est isomorphe \`a l'espace affine de dimension $\langle
\rho,\lambda+w\lambda\rangle$.
\end{lemme}

\dem
Pour $w=1$, l'assertion r\'esulte de mani\`ere \'evidente du lemme
pr\'ec\'edent, car on a les inclusions
$$J^\lambda
e_\lambda\subset S_\lambda\cap{\bar\Gr}_\lambda \subset
\vp^{\lambda}L^{<0}Ge_0\cap{\bar\Gr}_\lambda.$$

Pour $w\in W$ quelconque, on peut raisonner comme suit. D'apr\`es le lemme
pr\'ec\'edent, le morphisme
$$wJ^\lambda w^{-1}\rta \vp^{w\lambda}L^{<0}Ge_0\cap{\bar\Gr}_\lambda$$
d\'efini par $j\mapsto je_{w\lambda}$ est un isomorphisme.
Par ailleurs, la multiplication
$$(wJ^\lambda w^{-1}\cap LU)\times(wJ^\lambda w^{-1}\cap LU^-)\rta
wJ^\lambda w^{-1}$$
d\'efinit aussi un isomorphisme
si bien que pour $x\in L^{<w\lambda}U$ tel que
$x\vp^{w\lambda}\in{\bar\Gr}_\lambda$,
$x$ doit s'\'ecrire uniquement sous la forme $x=x_+ x_-$ avec
$x_+ \in wJ^\lambda w^{-1} \cap LU$
et $x_-\in wJ^\lambda w^{-1}\cap LU^-$.
Or, l'intersection $LU\cap LU^{-}$ est r\'eduite \`a l'\'el\'ement neutre
de sorte que
la seule possibilit\'e est $x=x_+$ et $x_-=1$.
Ceci prouve la premi\`ere assertion.

De plus, la multiplication induit un isomorphisme
$$
\prod_{\alpha\in R_+\cap w^{-1}R_+}
\prod_{i=0}^{\langle\alpha,\lambda\rangle-1}U_{w\alpha,i}
\stackrel{\cong}{\longrightarrow}
wJ^\lambda w^{-1}\cap LU. \eqno (*)
$$
Compte-tenu de l'\'egalit\'e
$$
\sum_{\alpha\in R_+\cap w^{-1}R_+}\, \alpha = \rho+ w^{-1}\rho,
$$
on obtient la seconde assertion. \findem

\medskip
On peut d\'eduire de ce lemme l'\'enonc\'e 3.1 dans le cas $\nu=w\lambda$
ainsi que 3.2 dans le cas $\nu=\lambda$.
En effet, l'inclusion \'evidente
$wJ^\lambda w^{-1}\cap LU\subset L^{\geq 0}U$ implique que
$S_{w\lambda}\cap{\bar\Gr}_\lambda$ est contenue dans l'orbite ouverte
$\Gr_\lambda$.
La restriction de $\calA_\lambda$ \`a
$S_{w\lambda}\cap{\bar\Gr}_\lambda$
est donc \'egale \`a :
$$
\calA_\lambda|_{S_{w\lambda}\cap{\bar\Gr}_\lambda}=
\Ql[\langle\rho,2\lambda \rangle]
(\langle\rho,\lambda \rangle).
$$
L'\'enonc\'e 3.1 dans le cas $\nu=w\lambda$
s'ensuit.

L'inclusion $J^\lambda\subset L^{\geq 0}U$ implique par
ailleurs que la restriction de $h$ \`a $J^\lambda$ est nulle.
L'\'enonc\'e 3.2 dans le cas $\nu=\lambda$ s'ensuit
donc aussi.

\smallskip
On aura besoin plus loin de l'\'enonc\'e plus
g\'en\'eral ci-dessous.
Pour tout $\sigma\in X^\smv_+$, on note $h_\sigma:LU\rta{\mathbb G}_a$
le morphisme $h_\sigma(x)=h(\vp^\sigma x\vp^{-\sigma})$ et aussi le
morphisme $h_\sigma:S_\lambda\rta{\mathbb G}_a$ qui s'en d\'eduit.
Du fait que $\sigma$ est dominant, la restriction de $h_\sigma$
\`a $L^{\geq 0}U$, et a fortiori \`a $J^\lambda$, est nulle.
Le lemme suivant r\'esulte donc \'egalement de la
discussion pr\'ec\'edente.

\begin{lemme}
Pour tous $\lambda, \sigma\in X^\smv_+$,  on a
$$\rmR\Gamma_c(S_\lambda,\calA_\lambda\otimes h_\sigma^*\calL_\psi)
=\Ql[-2\langle\rho,\lambda\rangle]
(-\langle\rho,\lambda\rangle).$$
\end{lemme}

\section{Minuscules}

Nous utilisons les notations fix\'ees dans la section 1.
Soit $\mu$ un \'el\'ement mi\-ni\-mal, non nul, de $X^\smv_+$. D'apr\`es le
lemme 1.1,
$\mu$ est un copoids {\it minuscule}, et l'on a l'\'enonc\'e suivant, que
nous rappelons
ici pour la commodit\'e du lecteur.

\begin{lemme} Soit $\mu$ minuscule.
On a $\Omega(\mu) = W\mu$.
Pour tout $\alpha\in R$, on a $\langle\alpha,\mu\rangle\in\{0,\pm 1\}$.
\end{lemme}

Si $\mu$ est minuscule, sa minimalit\'e
implique que l'orbite $\Gr_\mu$ est ferm\'ee.
Puisque tout \'el\'ement $\nu$ de $\Omega(\mu)$ est conjugu\'e \`a $\mu$
par l'action de $W$
les \'enonc\'es 3.1 et 3.2 sont
donc v\'erifi\'es pour $\lambda=\mu$ et $\nu\in\Omega(\mu)$.

On peut d\'ecrire explicitement la vari\'et\'e $\Gr_\mu$ et les strates
$S_{w\mu}\cap \Gr_\mu$.
On notera
$P$ le sous-groupe  de
$G$ engendr\'e par $T$ et les $U_\alpha$ avec
$\langle\alpha,\mu\rangle\leq 0$; c'est un sous-groupe parabolique
contenant $B^-$.

\begin{lemme}
On a un isomorphisme canonique $\Gr_\mu\rta G/P$ via lequel
$S_{w\mu}\cap\Gr_\mu$ s'identifie \`a $UwP/P$. \end{lemme}

\dem
Compte tenu du lemme 2.3 et de
la deuxi\`eme assertion du lemme 6.1, on sait que $L^{\geq 0}G\cap L^{\geq
\mu}G$ est
l'image inverse de $P$ par l'homo\-mor\-phisme $L^{\geq 0}G\rta G$.
On en d\'eduit l'isom\-orphisme
$$\Gr_\mu=L^{\geq 0}G/(L^{\geq 0}G\cap L^{\geq \mu}G)\cong G/P.$$

Compte tenu, de nouveau,
de la deuxi\`eme assertion du lemme 6.1, on sait que $J^\mu$ est \'egal \`a
$U_\mu^+=\prod_{\langle\alpha,\mu\rangle=1}U_\alpha$ qui est le radical
unipotent du sous-groupe parabolique
oppos\'e \`a $P$, et par cons\'equent
$$wJ^{\mu}w^{-1}\cap LU= wU_\mu^+ w^{-1}\cap U.$$ La seconde assertion du
lemme se d\'eduit alors du lemme 5.2. \findem

\section{Quasi-minuscules : \'etude g\'eom\'etrique}

Soit $\mu$ un copoids quasi-minuscule, c.\`a.d. un \'el\'ement minimal de
$X^\smv_+\setminus\{0\}$
minor\'e par $0$. Rappelons que, d'apr\`es le lemme 1.1, on a l'\'enonc\'e
suivant :

\begin{lemme}
Soit $\mu$ quasi-minuscule. Alors $\mu$ est \'egal \`a la coracine
$\gamma^\smv$ associ\'ee \`a une racine maximale $\gamma$.
On a $\Omega(\mu)=W\mu\cup\{0\}$. Pour toute racine $\alpha\in
R\setminus\{\pm\gamma\}$, on a
$\langle\alpha,\mu\rangle\in\{0,\pm 1\}$.
\end{lemme}

Puisque $0$ est le seul cocaract\`ere dominant qui minore $\mu$,
${\bar\Gr}_\mu$ est la r\'eunion de $\Gr_\mu$ et du point base
$e_0$.

D\'esignons encore par $P$ le sous-groupe parabolique de $G$ engendr\'e par
$T$ et par les sous-groupes radiciels $U_\alpha$ tels que
$\langle\alpha,\gamma^\vee\rangle\leq 0$.
Notons
$$
V={\mathfrak h}\oplus\bigoplus_{\alpha\in R\setminus\{\gamma\}}{\mathfrak
g}_\alpha
$$
o\`u $\mathfrak h$ est l'alg\`ebre de Lie de $T$ et o\`u ${\mathfrak g}_\alpha$
est le sous-espace de poids $\alpha$ de $\mathfrak g$.
D'apr\`es le lemme pr\'ec\'edent, $V$ est la somme des espaces de poids
$\nu$ dans
$\mathfrak g$ tels que $\langle\gamma, \nu\rangle \leq 1$. Il r\'esulte
alors de la d\'efinition
de $P$ que $V$ est $P$-stable.

Identifions ${\mathfrak g}_\gamma$ au quotient
${\mathfrak g}/V$ muni de sa structure de $P$-module. Consid\'erons le
fibr\'e en droites
$${\LL}_\gamma=G\times^{P}{\mathfrak g}_\gamma$$ au-dessus de $G/P$.

\begin{lemme}
L'orbite $\Gr_\mu$ est canoniquement isomorphe \`a ${\LL}_\gamma$.
\end{lemme}

\dem Par d\'efinition, le foncteur $R\mapsto G(R[\vp]/(\vp^2))$ est
repr\'esent\'e par le fibr\'e tangent $TG$ de $G$
qui est isomorphe au produit semi-direct $G\ltimes{\mathfrak g}$.
Compte tenu du lemme 2.3 et de la derni\`ere assertion du lemme 7.1, on sait que
$L^{\geq 0}G\cap L^{\geq\mu}G$
est exactement l'image inverse
de $P\ltimes V$ par l'homomorphisme canonique $L^{\geq 0}G\rta
G\ltimes{\mathfrak g}$.

Il en r\'esulte l'isomorphisme
$$\Gr_\mu\cong (G\ltimes{\mathfrak g})/(P\ltimes V)
=G\times^P({\mathfrak g}/V)$$
dont le terme de droite n'est autre que ${\LL}_\gamma$. \findem

Le fibr\'e $\LL_\gamma$ se compactifie de fa{\c c}on naturelle en
un fibr\'e en droites projectives. En effet, on a $$\LL_\gamma
\hookrightarrow\Proj(\LL_\gamma \oplus \calO_{G/P}) =\PP_\gamma.$$
Soit $\LL_{-\gamma} = G\times^P {\mathfrak g}_{-\gamma}$ le fibr\'e dual.
On a un isomorphisme naturel
$$\Proj(\LL_\gamma \oplus \calO_{G/P}) \cong \Proj(\calO_{G/P} \oplus
\LL_{-\gamma})
=\PP_{-\gamma},$$
si bien qu'on peut voir $\PP_\gamma$ comme la r\'eunion de $\LL_\gamma$ et
de $\LL_{-\gamma}$.
Notons $\phi_{\pm\gamma}$ le morphisme $\LL_{\pm\gamma}\to G/P$ et
$\epsilon_{\pm\gamma}$ la section nulle $G/P\to \LL_{\pm\gamma}$. On a
alors
$$\PP_\gamma = \LL_{\pm\gamma}\cup \epsilon_{\mp\gamma}(G/P).$$

\begin{lemme}
L'isomorphisme $\LL_\gamma\rta\Gr_\mu$ se prolonge en un morphisme
$$\pi_\gamma:\PP_\gamma\rta{\bar\Gr}_\mu$$ qui envoie
$\epsilon_{-\gamma}(G/P)$ sur le point $e_0$.
\end{lemme}

\dem
Il s'agit d'un cas particulier du th\'eor\`eme principal de Zariski. Nous
reproduisons l'argument classique pour la commodit\'e du lecteur. Notons
$\Gamma$ l'adh\'erence du graphe de l'isomorphisme $\LL_\gamma\rta\Gr_\mu$
dans $\PP_\gamma\times{\bar\Gr}_\mu$. Puisque le compl\'ementaire de
${\Gr}_\mu$ dans ${\bar\Gr}_\mu$ est constitu\'e d'un seul point $e_0$,
toutes les fibres de la projection
$\Gamma\rta \PP_\gamma$ ne contiennent qu'un point. En particulier, cette
projection est quasi-finie. Puisque de plus, elle est propre, elle est
finie. Mais $\PP_\gamma$ est lisse, en particulier normale, donc le morphisme
fini, birationnel $\Gamma\rta\PP_\gamma$ doit \^etre un isomorphisme. En
inversant cet isomorphisme et en le composant avec l'autre projection
$\Gamma\rta{\bar\Gr}_\mu$, on obtient
le morphisme $\pi_\gamma:\PP_\gamma\rta{\bar\Gr}_\mu$ voulu. \findem

Signalons que l'action de $L^{\geq 0}G$ sur $\Gr_\mu$ se prolonge
\`a ${\mathbb P}_\gamma$ et que la r\'esolution
$\pi_\gamma$ est \'equivariante par rapport \`a cette action.
Nous n'utiliserons cette information que dans la remarque
situ\'ee apr\`es le corollaire 9.7, laquelle ne sert pas
dans le reste de l'article.

\bigskip
On a la description explicite suivante des strates
$S_{w\mu}\cap{\bar\Gr}_\mu$.

\begin{lemme}
Si $w\gamma\in R_+$ alors
$$S_{w\mu}\cap{\bar\Gr}_\mu=\phi_\gamma^{-1}(UwP/P).$$ Si $w\gamma\in R_-$ alors
$$S_{w\mu}\cap{\bar\Gr}_\mu=\epsilon_\gamma(UwP/P).$$ \end{lemme}

\dem
D'apr\`es la formule $(*)$ \'etablie dans la d\'emonstration
du lemme 5.2, l'on a:
$$
wJ^\mu w^{-1}\cap LU =
\prod_{\alpha\in R_+\cap w^{-1}R_+}
\prod_{i=0}^{\langle\alpha,\mu\rangle-1}U_{w\alpha,i}.
$$
De plus, comme $\langle\alpha, \mu\rangle \leq 1$ pour tout $\alpha\in
R_+\setminus\{\gamma\}$,
d'apr\`es le lemme 7.1, on obtient que ce produit est \'egal \`a
$$
U_{w\gamma,1}
\prod_{\alpha\in R_+\cap w^{-1}R_+}U_{w\alpha,0}
$$
si $w\gamma\in R_+$, et \`a
$$
\prod_{\alpha\in R_+\cap w^{-1}R_+}U_{w\alpha,0}
$$
sinon.
Le lemme s'en d\'eduit.
\findem

On note $W_\gamma$ le stabilisateur de $\gamma$ dans $W$, et
$\Delta_\gamma$ l'ensemble des racines simples conjugu\'ees \`a $\gamma$.

\begin{corollaire} On a une stratification
$$
S_0\cap{\bar\Gr}_\mu = \{e_0\}\cup
\bigcup_{\deuxind{ w\in W/W_\gamma}{w\gamma\in R_-}}
\phi_\gamma^{-1}(UwP/P)\setminus \epsilon_\gamma(UwP/P).
$$
En particulier, les composantes irr\'eductibles de $S_0\cap{\bar\Gr}_\mu$
sont en
bijection avec $\Delta_\gamma$ et sont toutes de dimension
$\langle\rho,\mu\rangle$.

On a aussi la stratification
$$\pi^{-1}_\gamma(S_0\cap{\bar\Gr}_\mu)
=\bigcup_{\deuxind{ w\in W/W_\gamma}{ w\gamma\in R_-}}
\phi_{-\gamma}^{-1}(UwP/P)\cup
\bigcup_{\deuxind{w\in W/W_\gamma}{w\gamma\in R_+}}
\epsilon_{-\gamma}(UwP/P).$$
\end{corollaire}

\section{Quasi-minuscules : \'etude cohomologique}

Les notations de la section pr\'ec\'edente restent en vigueur. En particulier,
$\mu=\gamma^{\smv}$ est quasi-minuscule. La r\'esolution
$\pi_\gamma:\PP_\gamma\rta{\bar\Gr}_\mu$ permet de calculer la
cohomologie d'intersection locale de $\calA_\mu$ en le point singulier
isol\'e $e_0$.
L'\'enonc\'e suivant est d\^u \`a Kazhdan et Lusztig \cite[Lemme
4.5]{Kazhdan-Lusztig}.
En fait, dans notre situation les hypoth\`eses sont un peu plus faibles, mais
leur argument s'applique encore. Nous d\'etaillons la d\'emonstration
pour la commodit\'e du lecteur.

\begin{lemme}
Soit $d=2\langle\rho,\mu\rangle$ la dimension de ${\bar\Gr}_\mu$.
Pour $i\geq 0$, le groupe $\rmH^{i}(\calA_\mu)_{e_0}$ est nul. Pour $i< 0$,
on a une suite exacte courte
$$
0\rta
\rmH^{i+d-2}(G/P)(d/2-1)\,\hfld{\wedge\,c_{-\gamma}}{}\,
\rmH^{i+d}(G/P)(d/2)\rta \rmH^{i}(\calA_\mu)_{e_0}\rta 0,
$$
o\`u
$c_{-\gamma}\in H^2(X_\gamma)(1)$ est la classe de Chern de ${\mathbb
L}_{-\gamma}$.
\end{lemme}

\dem
Notons ${\bar\Gr}'_\mu$ l'ouvert de ${\bar\Gr}_\mu$
$${\bar\Gr}'_\mu={\bar\Gr}_\mu-\pi_\gamma\circ\epsilon_\gamma(G/P) ;$$
on a $\pi_\gamma^{-1}({\bar\Gr}'_\mu)=\LL_{-\gamma}$. Notons $\calA'_\mu$
la restriction de $\calA_\mu$ \`a cet ouvert. Notons
$\imath$ l'inclusion du point ferm\'e $\imath:\{e_0\}\rta {\bar\Gr}'_\mu$.
La fl\`eche naturelle $\calA'_\mu\rta
\imath_*\imath^*\calA'_\mu$ induit le morphisme de restriction sur la
cohomologie (sans support)
$$\imath^*:\rmR\Gamma({\bar\Gr}'_\mu,\calA'_\mu)\rta (\calA'_\mu)_{e_0}.$$
On d\'emontre d'abord que
$\imath^*$ est un isomorphisme.

Pour cela, utilisons le th\'eor\`eme de d\'ecomposition de
Beilinson, Bernstein, Deligne et Gabber \cite{BBD}.
Puisque $\pi_\gamma:\PP_\gamma\rta{\bar\Gr}_\mu$ est un isomorphisme en
dehors de
$e_0$, on a une d\'ecomposition
$$\rmR\pi_{\gamma,*}\Ql[d](d/2)=\calA_\mu\oplus\calC,$$
o\`u $\calC$ est un complexe support\'e par le point $e_0$.

La section nulle $\epsilon_{-\gamma}:G/P\rta\LL_{-\gamma}$ d\'efinit le
morphisme de restriction
$$
\rmR\Gamma(\LL_{-\gamma},\Ql)\rta\rmR\Gamma(G/P,\Ql)
$$
qui est un isomorphisme
puisque $\LL_{-\gamma}$ est
un fibr\'e en droites au-dessus de $G/P$. Or, ce morphisme
est la somme directe du morphisme identit\'e $\calC\rta\calC$ avec le
morphisme $\imath^*:\rmR\Gamma({\bar\Gr}'_\mu,\calA'_\mu)\rta
(\calA'_\mu)_{e_0}$.
Ce dernier est donc lui aussi un isomorphisme.

Pour $i\geq 0$, l'annulation $\rmH^i(\calA_\mu)_{e_0}=0$ fait partie des
propri\'et\'es qui carac\-t\'erisent
le complexe d'intersection $\calA_\mu$.
Cette annulation implique,
via la suite exacte longue de cohomologie \`a support, que la fl\`eche
$$\rmH^{i+d}_c(\LL_\gamma^\times)(d/2)\rta \rmH^{i}_c(\Gr'_\mu,\calA'_\mu)$$
est un isomorphisme d\`es que $i>0$. Par dualit\'e de Poincar\'e, la fl\`eche
$$\rmH^{i}(\Gr'_\mu,\calA'_\mu)\rta
\rmH^{i+d}(\LL_\gamma^\times)(d/2)$$
est un isomorphisme pour $i<0$.
Dans la suite exacte longue de Wang
$$\rta\rmH^{i+d-2}(G/P)(i/2-1)\hfld{\wedge\,c_{-\gamma}}{}
\rmH^{i+d}(G/P)(i/2)\rta
\rmH^{i+d}(\LL_\gamma^\times)(i/2)\rta$$
la fl\`eche $\rmH^{i+d-2}(G/P)\hfld{\wedge\,c_{-\gamma}}{}
\rmH^{i+d}(G/P)$ est injective pour $i\leq 0$
d'apr\`es le th\'eor\`eme
de Lefschetz   difficile \cite{Weil II}. On en d\'eduit les suites exactes
courtes
$$0\rta \rmH^{i+d-2}(G/P)(-1)\hfld{\wedge\,c_{-\gamma}}{}
\rmH^{i+d}(G/P)\rta \rmH^{i+d}(\LL_\gamma^\times)\rta 0$$ pour $i<0$.
Le lemme est d\'emontr\'e. \findem

\begin{corollaire}
Soit $\calC$ le facteur support\'e par $e_0$ dans la d\'ecomposition
$$\rmR\pi_{\gamma,*}\Ql[d](d/2)=\calA_\mu\oplus\calC.$$
Pour $i<0$, on a
$$\rmH^i(\calC)=\rmH^{i+d-2}(G/P)(d/2-1).$$ Pour $i \geq 0$, on a
$$\rmH^i(\calC)=\rmH^{i+d}(G/P)(d/2).$$
\end{corollaire}

On peut maintenant d\'emontrer l'\'enonc\'e 3.1 dans le cas o\`u $\lambda$
est un cocaract\`ere quasi-minuscule
$\mu=\gamma^{\smv}$. Compte tenu de la discussion
qui suit le lemme 5.2, il ne reste plus qu'\`a traiter le cas $\nu=0$.

\begin{lemme}
On a un isomorphisme
$$\rmR\Gamma_c(S_0,\calA_\mu) \cong \Ql^{|\Delta_\gamma|}.$$
\end{lemme}

\dem
D'apr\`es le th\'eor\`eme de changement de base pour un morphisme propre,
on a
$$
\rmR\Gamma_c(\pi_\gamma^{-1}(S_0\cap{\bar\Gr}_{\mu}),\Ql)[d](d/2)
=\rmR\Gamma_c(S_0,\calA_\mu)\oplus{\mathcal C}.
$$
Rappelons la stratification
obtenue en 7.5
$$
\pi_\gamma^{-1}(S_0\cap{\bar\Gr}_\mu)=
\bigcup_{\deuxind{w\in W/W_\gamma}{w\gamma\in R_-}}
\phi^{-1}_{-\gamma}(UwP/P)\cup
\bigcup_{\deuxind{w\in W/W_\gamma}{ w\gamma\in R_+}}
\epsilon_{-\gamma}(UwP/P).
$$

D'apr\`es le lemme 5.2, chaque strate $S_{w\mu} \cap {\bar\Gr}_{\mu}$ est
de dimension
$$
\langle\rho,w\mu+\mu\rangle.
$$
De plus, d'apr\`es le lemme 7.4, l'on a
$$
S_{w\mu} \cap {\bar\Gr}_{\mu} = \begin{cases} \phi_\gamma^{-1}(UwP/P)
\text{ si $w\gamma\in R_+$}; \\
\phantom{O}\epsilon_\gamma(UwP/P) \text{ si $w\gamma\in R_-$}.\end{cases}
$$
On en d\'eduit que
si $w\gamma\in R_-$, alors la strate $\phi^{-1}_{-\gamma}(UwP/P)$
est un espace affine de dimension
$$
\langle\rho,w\mu+\mu\rangle + 1.
$$
On a, dans ce cas,
l'in\'egalit\'e
$\dim(\phi^{-1}_{-\gamma}(UwP/P))\leq d/2$ qui devient une \'egalit\'e si
et seulement si $w\gamma$ est l'oppos\'e d'une racine simple.

D'autre part, si $w\gamma\in R_+$, alors
la strate $\epsilon_{-\gamma}(UwP/P)$ est un espace affine de dimension
$$
\langle\rho,w\mu+\mu\rangle - 1.
$$
On a, dans ce cas, l'in\'egalit\'e
$\dim(\epsilon_{-\gamma}(UwP/P))\geq d/2$ qui devient une \'egalit\'e si et
seulement si $w\gamma$ est une
racine simple.

Pour $i<0$, on a donc
$$\dim \rmH_c^{i+d}(\pi_\gamma^{-1}(S_0\cap{\bar\Gr}_{\mu}))
=\dim \rmH^i(\calC)$$
et ces deux nombres valent
$$|\{w\in W/W_\gamma\mid\langle\rho,w\mu+\mu\rangle=(i+d)/2-1\}|$$ Pour
$i>0$, on a aussi
$$\dim \rmH_c^{i+d}(\pi_\gamma^{-1}(S_0\cap{\bar\Gr}_{\mu}))
=\dim \rmH^i(\calC)$$
et ces deux nombres valent
$$|\{w\in W/W_\gamma\mid\langle\rho,w\mu+\mu\rangle=(i+d)/2+1\}|.$$

Pour $i=0$, on a
$$\dim \rmH^d_c(\pi_\gamma^{-1}(S_0\cap{\bar\Gr}_{\mu})) =2\, |\Delta_\gamma|$$
et $$\dim(\rmH^0(\calC))=|\Delta_\gamma|.$$ Le lemme s'en d\'eduit. \findem

D\'emontrons maintenant l'\'enonc\'e 3.2 dans le cas $\nu=0$ et
$\lambda=\mu$ quasi-minuscule. Nous d\'emontrons en fait un \'enonc\'e un
peu plus g\'en\'eral. Rappelons que pour tout $\sigma\in X^\smv$,
on a d\'efini un morphisme $h_\sigma:S_0\rta{\mathbb G}_a$, voir 5.3.

\begin{lemme}
Pour tout $\sigma\in X^\smv_+$, on a un isomorphisme
 $$\rmR\Gamma_c(S_0,
\calA_{\mu}\!\otimes h_\sigma^{*}\calL_\psi) =\Ql^{|\Delta_\gamma^\sigma|},$$
o\`u $\Delta_\gamma^\sigma$ est l'ensemble des
$\alpha\in\Delta_\gamma$ telles que $\langle\alpha,\sigma\rangle >0$.
\end{lemme}

La d\'emonstration de 8.4 suit le m\^eme sch\'ema que celle du lemme 8.3
qui est, par ailleurs, un cas particulier de 8.4. Il suffit, en fait, de
d\'emontrer l'\'enonc\'e g\'eom\'etrique
suivant.

\begin{lemme}
\begin{enumerate}
\item
Les restrictions de $h_\sigma\circ\pi_\gamma$ aux strates
$\epsilon_{-\gamma}(UwP/P)$ avec $w\gamma\in R_+$ sont nulles.
\item
Les restrictions aux strates $\phi^{-1}_{-\gamma}(UwP/P)$ avec $w\gamma\in R_-$
sont nulles aussi, \`a l'exception des $\phi^{-1}_{-\gamma}(UwP/P)$ telles
que $-w\gamma$ est une
racine simple orthogonale \`a $\sigma$.
\item
Les restrictions \`a ces
derni\`eres sont non nulles et lin\'eaires par rapport
\`a la structure restreinte du fibr\'e en droites $\LL_{-\gamma}$.
\end{enumerate}
\end{lemme}

\dem
La premi\`ere assertion est \'evidente parce que toutes les strates
$\epsilon_{-\gamma}(UwP/P)$ sont contenues dans $\pi_\gamma^{-1}(e_0)$.
Pour les deux derni\`eres assertions, il suffit naturellement de calculer
les restrictions de $h_\sigma$ \`a
$$\phi^{-1}_{-\gamma}(UwP/P)\setminus\epsilon_{-\gamma}(UwP/P)=
\phi^{-1}_{\gamma}(UwP/P)\setminus\epsilon_{\gamma}(UwP/P).$$
Pour tout $w\in W$,
on a un isomorphisme $$\phi^{-1}_{\gamma}(UwP/P)\setminus
\epsilon_{\gamma}(UwP/P)
=UwP/P\times{\mathbb G}_m$$
puisque tout point
$y\in\phi^{-1}_{\gamma}(UwP/P)\setminus\epsilon_{\gamma}(UwP/P)$
s'\'ecrit d'une fa{\c c}on unique sous la forme
$$y=uwU_{\gamma,1}(x)\vp^\mu e_0$$
avec $x\in{\mathbb G}_m$ et $u\in U\cap w^{-1}U^+_\gamma w$, o\`u
$U^+_\gamma$ est le radical unipotent du parabolique (positif) oppos\'e \`a
$P$.

En faisant passer $w$ vers la droite, on obtient $$uwU_{\gamma,1}(x)\vp^\mu
e_0=uU_{w\gamma,1}(x)\vp^{w\mu}e_0.$$
Posons $t=-\vp x$ et $\alpha=w\gamma$ et r\'ecrivons le membre de droite
avec ces nouvelles notations
$$uU_{w\gamma}(\vp x)\vp^{w\mu}e_0=u U_\alpha(-t) t^{\alpha^\smv} e_0.$$

Rappelons la relation de Steinberg \cite[chap. 3, lemme 19]{Steinberg}
$$t^{\alpha^\smv}w_\alpha=U_\alpha(t)U_{-\alpha}(-t^{-1})U_\alpha(t)$$ qui
vaut pour toute racine $\alpha$, pour tout $t$
inversible, et o\`u le repr\'esentant $w_\alpha\in G$ est ind\'ependant de
$t$. On a donc
$$uU_\alpha(-t) t^{\alpha^\smv} e_0=U_{-\alpha}(-t^{-1})U_\alpha(t)
w_\alpha^{-1}e_0.$$
Or $U_\alpha(t)w_\alpha^{-1}e_0=e_0$ si bien que $$uU_\alpha(-t)
t^{\alpha^\smv} e_0=uU_{-\alpha,-1}(x^{-1})e_0.$$

Puisque $u\in U$, on a $h(u)=0$ de sorte que
$$h(uU_{-\alpha,-1}(x^{-1}))=h(U_{-\alpha,-1}(x^{-1})).$$ Si $-\alpha$
n'est pas une racine simple, on a $h_\sigma(U_{-\alpha,-1}(x^{-1}))=0$.
Si $-\alpha\in\Delta$ mais $\langle-\alpha,\sigma\rangle>0$, on a aussi
$h_\sigma(U_{-\alpha,-1}(x^{-1}))=0$. Dans le cas o\`u $-\alpha\in\Delta$
est une racine simple orthogonale \`a $\sigma$, on a
$h_\sigma(U_{-\alpha,-1}(x^{-1}))=x^{-1}$. \findem

\section{Convolution}
Rappelons que $M$ est l'ensemble des \'el\'ements minimaux dans
$X^\smv_+\setminus\{0\}$.
Pour une suite $\mu_\bullet=(\mu_1,\ldots,\mu_n)$ d'\'el\'ements
de $M$, on consid\`ere le sous-sch\'ema projectif $${\bar\Gr}_{\mu_\bullet}=
{\bar\Gr}_{\mu_1}\tildetimes\cdots\tildetimes{\bar\Gr}_{\mu_n}$$ de
$\Gr^n$.
La projection sur le dernier facteur de $\Gr^n$
d\'efinit un morphisme propre
$$
m_{\mu_\bullet} : {\bar\Gr}_{\mu_\bullet} \rta{\bar\Gr}_{|\mu_\bullet|},
$$
o\`u $|\mu_\bullet|=\mu_1+\cdots+\mu_n$.

Soit $\nu_\bullet=(\nu_1,\ldots,\nu_n)$ une suite d'\'el\'ements de $X^\smv$.
Pour tout $i=1,\ldots,n$, posons $\sigma_i=\nu_1+\cdots+\nu_i$. Notons
$S_{\nu_\bullet}\cap {\bar\Gr}_{\mu_\bullet}$ l'intersection
$$S_{\nu_\bullet}\cap {\bar\Gr}_{\mu_\bullet}
=S_{\sigma_1}\times\cdots\times S_{\sigma_n} \cap {\bar\Gr}_{\mu_\bullet}$$
dans $\Gr^n$. Il est clair que les
$S_{\nu_\bullet}\cap {\bar\Gr}_{\mu_\bullet}$ forment une stratification de
${\bar\Gr}_{\mu_\bullet}$.

\begin{lemme}
On a un isomorphisme canonique
$$S_{\nu_\bullet}\cap {\bar\Gr}_{\mu_\bullet}
\cong(S_{\nu_1}\cap{\bar\Gr}_{\mu_1})\times\cdots\times
(S_{\nu_n}\cap{\bar\Gr}_{\mu_n}).$$
\end{lemme}

\dem
On montre facilement par r\'ecurrence que tout point $$(y_1,\ldots,y_n)\in
S_{\nu_\bullet}\cap {\bar\Gr}_{\mu_\bullet}$$ s'\'ecrit de mani\`ere unique
sous la forme $$\begin{array}{rcl}
y_1 & = & x_1\vp^{\nu_1}e_0 \\
&\ldots & \\
y_n & = & x_1\vp^{\nu_1}\ldots x_n\vp^{\nu_n}e_0\\ \end{array}$$
avec $x_i\in L^{<\nu_i}U$ tel que
$x_i\vp^{\nu_i}e_0\in{\bar\Gr}_{\mu_i}$. Le lemme s'en d\'eduit.\findem

Pour que la strate $S_{\nu_\bullet}\cap {\bar\Gr}_{\mu_\bullet}$ soit non
vide, il est donc
n\'ecessaire que, pour tout
$i=1,\ldots,n$, $\nu_i$ appartienne \`a $\Omega(\mu_i)$.

\begin{corollaire}
Soient $\mu_1,\ldots,\mu_n$ des \'el\'ements de $M$.
Pour toute suite $\nu_\bullet$ avec $\nu_i\in\Omega(\mu_i)$, toutes les
composantes de $S_{\nu_\bullet}\cap {\bar\Gr}_{\mu_\bullet}$ sont de dimension
$\langle\rho,|\nu_\bullet|+|\mu_\bullet|\rangle$.
\end{corollaire}

\dem
D'apr\`es le lemme 5.2 et le corollaire 7.5, chaque
$S_{\nu_i}\cap{\bar\Gr}_{\mu_i}$
est purement de dimension $\langle\rho,\nu_i+\mu_i\rangle$.
Le corollaire d\'ecoule donc du lemme pr\'ec\'edent.
\findem

En fait, pour $\lambda\in X^\smv_+$ et $\nu\in\Omega(\lambda)$ arbitraires,
$S_{\nu}\cap{\bar\Gr}_{\mu}$ est
purement de dimension
$\langle\rho,\nu+\mu\rangle$. Ce r\'esultat est \'enonc\'e dans
\cite[4.5]{Mirkovic-Vilonen}
avec seulement quelques indications de d\'emonstration.
Nous avons pu en \'etablir une autre d\'emonstration, en utilisant
l'approche en termes de
repr\'esentations d'alg\`ebres de Lie affines.
Signalons aussi qu'on peut d\'eduire cette formule de dimension, sans
l'assertion de pure dimension, du th\'eor\`eme 3.1.

L'\'enonc\'e suivant est aussi un cas particulier d'un autre lemme \'enonc\'e
dans \cite[2.6]{Mirkovic-Vilonen}. Nous proposons ici une d\'emonstration
un peu diff\'erente.

\begin{lemme}
Pour tout $\lambda\in X^\smv_+$ avec $\lambda\leq |\mu_\bullet|$, on a
$$
\dim(m_{\mu_\bullet}^{-1}(\Gr_\lambda))\leq\langle\rho,\lambda+|\mu_\bullet|
\rangle.
$$
Autrement dit, le morphisme $m_{\mu_\bullet}$ est semi-petit. \end{lemme}

\dem
Puisque $S_\mu\cap\Gr_\mu$ est un ouvert dense de $\Gr_\mu$, d'apr\`es 5.1
et 5.2, et puisque les fibres de $m_{\mu_\bullet}$ sont toutes isomorphes
les unes aux
autres au-dessus de l'orbite $\Gr_\mu$, on a
$$
\dim(m^{-1}_{\mu_\bullet}(\Gr_\mu))
=\dim(m^{-1}_{\mu_\bullet}(S_\mu\cap\Gr_\mu)).
$$
Par cons\'equent, on a
$$
\dim(m^{-1}_{\mu_\bullet}(\Gr_\mu))
\leq\dim(m^{-1}_{\mu_\bullet}(S_\mu\cap{\bar\Gr}_\lambda)).
$$
Or, on a une stratification
$$
m^{-1}_{\mu_\bullet}(S_\mu\cap{\bar\Gr}_\lambda)
=\bigcup_{|\nu_\bullet|=\lambda}S_{\nu_\bullet}\cap
{\bar\Gr}_{\mu_\bullet},
$$
o\`u toutes les strates sont de dimension
$\langle\rho,\lambda+|\mu_\bullet|\rangle$, d'o\`u le lemme.
\hfill$\square$

\begin{proposition}
Le produit de convolution $\calA_{\mu_1}*\cdots*\calA_{\mu_n}$ est un
faisceau pervers. Il se d\'ecompose en somme directe
$$
\calA_{\mu_1}*\cdots*\calA_{\mu_n}=\bigoplus_{\lambda\leq|\mu_\bullet|}
\calA_\lambda\otimes V^{\lambda}_{\mu_\bullet},
$$
o\`u les
$V^{\lambda}_{\mu_\bullet}$ sont des $\Ql$-espaces vectoriels dont la
dimension vaut le nombre de composantes irr\'eductibles de
$m^{-1}_{\mu_\bullet}(S_\lambda\cap{\bar\Gr}_{\mu_\bullet})$ qui sont
enti\`erement
contenues dans $m^{-1}_{\mu_\bullet}(S_\lambda\cap{\bar\Gr}_\lambda)$.
\end{proposition}

\dem
Si les $\mu_i$ sont tous minuscules, le sch\'ema source est lisse.
Comme $m_{\mu_\bullet}$ est semi-petit,
l'image directe $\rmR (m_{\mu_\bullet}){}_*\Ql
[\dim(\bar\Gr_{\mu_\bullet})]$ est perverse.

En g\'en\'eral, on a la stratification du sch\'ema source
$$
{\bar\Gr}_{\mu_\bullet}=\bigcup_{\mu'_i\leq\mu_i}
\Gr_{\mu'_1}\tildetimes\cdots\tildetimes\Gr_{\mu'_n},
$$ o\`u, comme
$\mu_i\in M$, chaque $\mu'_i$ est ou bien \'egal \`a $\mu_i$,
ou bien \'egal \`a $0$. Dans tous les cas, le lemme pr\'ec\'edent
s'applique encore \`a ${\bar\Gr}_{\mu'_\bullet}$ et nous permet d'obtenir
l'in\'egalit\'e
$$
\dim(m^{-1}_{\mu_\bullet}(\Gr_\lambda)\cap {\bar\Gr}_{\mu'_\bullet})
\leq\langle\rho,\lambda+|\mu'_\bullet|\rangle$$ pour tout $\lambda\leq
|\mu'_\bullet|$.

Par ailleurs, ${\bar\Gr}_{\mu_\bullet}$ muni de la stratification par les
$\Gr_{\mu'_1}\tildetimes\cdots\tildetimes\Gr_{\mu'_n}$, est localement
isomorphe \`a ${\bar\Gr}_{\mu_1}\times\cdots \times{\bar\Gr}_{\mu_1}$ muni
de la stratification par les $\Gr_{\mu'_1}\times\cdots\times\Gr_{\mu'_n}$,
voir la section 2. Par cons\'equent, pour $\mu'_\bullet<\mu_\bullet$,
$$\rmH^i({\rm IC}({\bar\Gr}_{\mu_\bullet})
|_{\Gr_{\mu'_1}\tildetimes\cdots\tildetimes\Gr_{\mu'_n}})$$ s'annule
d\`es que $i\geq -2\langle\rho,|\mu'_\bullet|\rangle$.

Il r\'esulte de ces deux derni\`eres assertions que
$\rmR (m_{\mu_\bullet}){}_*{\rm IC}({\bar\Gr}_{\mu_\bullet})$
est un faisceau pervers. Sa d\'ecomposition
\cite{BBD}, a priori sur ${\bar k}$, doit avoir la forme
$$\calA_{\mu_1}*\cdots*\calA_{\mu_n}=\bigoplus_{\lambda\leq|\mu_\bullet|}
\calA_\lambda\otimes V_{\mu_\bullet}^\lambda$$
o\`u les $V^{\lambda}_{\mu_\bullet}$ sont des espaces vectoriels
de dimension finie sur $\Ql$,
parce que tous ses facteurs directs sont aussi $L^{\geq 0}G$-\'equivariants.

Les espaces vectoriels $V^{\lambda}_{\mu_\bullet}$ admettent une base
canonique index\'ee par les composantes irr\'eductibles de
$m^{-1}_{\mu_\bullet}(\Gr_\lambda)$ de dimension exactement
$\langle\rho,|\mu_\bullet|-\lambda\rangle$. D'apr\`es la d\'emonstration du
lemme 9.3, celles-ci correspondent bijectivement aux composantes
irr\'eductibles
de $m^{-1}_{\mu_\bullet}(S_\lambda\cap{\bar\Gr}_{\mu_\bullet})$ qui sont
enti\`erement
contenues dans $m_{\mu_\bullet}^{-1}(S_\lambda\cap{\bar\Gr}_\lambda)$.

Puisque ces composantes sont toutes d\'efinies sur $k$, la d\'ecomposition
est en fait valable sur $k$. \findem

Soit $\mu_\bullet$ comme pr\'ec\'edemment, c.\`a.d. $\mu_\bullet$
est une
suite $(\mu_1,\dots,\mu_n)$ d'\'el\'ements de $M$. A la suite de Littelmann
\cite{Littelmann}, on
appellera
$\mu_\bullet$-{\it chemin} une donn\'ee combinatoire $\chi$ du type suivant
\begin{itemize}
\item une suite de sommets $\sigma_1,\ldots,\sigma_n$ dans $X^\smv$ tels
que pour tout $i=1,\ldots,n$, on a
$\nu_i=\sigma_i-\sigma_{i-1}\in\Omega(\mu_i)$ ; \item des applications
$$p_i:[0,1]\rta X^\smv\otimes_{\mathbb Z}{\mathbb R}$$ v\'erifiant
\begin{itemize}
\item si $\sigma_{i-1}\not=\sigma_i$, alors
$$p_i(t)=(1-t)\sigma_{i-1}+t\sigma_i$$
\item si $\sigma_{i-1}=\sigma_i$, alors
$$
p_i(t)=\begin{cases}\sigma_{i-1}-t\alpha^\smv_i &{\rm pour\ \ } 0\leq t\leq
1/2, \\
\sigma_{i-1}+(t-1)\alpha^\smv_i & {\rm pour\ \ } 1/2\leq t\leq 1,
\end{cases}
$$
o\`u $\alpha^\smv_i$ est une coracine simple conjugu\'ee \`a $\mu_i$.
\end{itemize}
\end{itemize}

En mettant bout \`a bout les images des $p_i$, on obtient un chemin dans
$X^\smv\otimes_{\mathbb Z}{\mathbb R}$ allant de $0$ \`a $\sigma_n$.
Le $\mu_\bullet$-chemin $\chi$ est dit {\it dominant} s'il est enti\`erement
contenu dans la chambre dominante $(X^\smv\otimes_{\mathbb Z}{\mathbb R})_+$.

D'apr\`es le lemme 5.2, chaque $S_{w\mu_i}\cap{\bar\Gr}_{\mu_i}$ est
irr\'eductible.
De plus, d'apr\`es le corollaire 7.5,
si $\mu_i$ est quasi-minuscule, disons $\mu_i = \gamma^\smv_i$,  et si
$\nu=0$, alors l'ensemble des composantes irr\'eductibles de
$S_0\cap{\bar\Gr}_{\mu_i}$ est en bijection canonique avec l'ensemble
$\Delta_\gamma$ des racines simples $\alpha$ conjugu\'ees \`a
$\gamma$. Compte tenu du lemme 9.1, pour tout $\nu\in\Omega(|\mu_\bullet|)$,
l'ensemble des composantes irr\'eductibles de
$\pi^{-1}(S_\nu\cap{\bar\Gr}_{|\mu_\bullet|})$ est en bijection canonique
avec l'ensemble des $\mu_\bullet$-chemins $\chi$ allant de $0$ \`a $\nu$.
Notons $C_\chi$ la composante correspondante \`a $\chi$.

\begin{lemme}
Soient $\nu\in\Omega(|\mu_\bullet|)$ dominant et $\chi$ un
$\mu_\bullet$-chemin dominant allant de $0$ \`a $\nu$.
Alors la composante
$C_\chi$ est contenue dans $\pi^{-1}(S_\nu\cap{\bar\Gr}_\nu)$.
\end{lemme}

\dem
Notons $I(\chi)$ l'ensemble des indices $i=1,\ldots,n$ tels que
$\sigma_{i-1}=\sigma_i$.

Si $i\notin I(\chi)$, $\nu_i$ est non nul et est donc conjugu\'e \`a $\mu_i$.
D'apr\`es le lemme 5.2, un point
$p_i\in S_{\nu_i}\cap{\bar\Gr}_{\mu_i}$
s'\'ecrit de mani\`ere unique sous la forme $p_i=u_i\vp^{\nu_i}e_0$ avec
$u_i\in wJ^{\mu_i}w^{-1}\cap LU$. En particulier, $u_i\in L^{\geq 0}U$.

Si $i\in I(\chi)$, alors $\mu_i$ est quasi-minuscule, disons $\mu_i =
\gamma^\smv_i$, et
l'hypo\-th\`e\-se $\chi$ dominant implique
$\langle\alpha_i,\sigma_{i-1}\rangle\geq 1$. D'apr\`es le corollaire 7.5,
la composante irr\'eductible de $S_0\cap {\bar\Gr}_{\gamma_i^\smv}$
correspondant \`a $\alpha_i=w_i\gamma_i$, contient comme ouvert dense le
${\mathbb G}_m$-torseur trivial
$$\phi_{\gamma_i}^{-1}(Uw_iP_i/P_i)
\setminus\epsilon_{\gamma_i}^{-1}(Uw_iP_i/P_i)$$
au-dessus de $Uw_iP_i/P_i$. D'apr\`es la d\'emonstration du lemme 8.5,
pour $i\in I(\chi)$ chaque point
$$p_i\in \phi_{\gamma_i}^{-1}(Uw_iP_i/P_i)
\setminus\epsilon_{\gamma_i}^{-1}(Uw_iP_i/P_i)$$
s'\'ecrit de mani\`ere unique sous
la forme $uU_{\alpha_i,-1}(x)e_0$ avec $u\in U\cap w^{-1}U^+_{\gamma_i}w$
et $x\in{\mathbb G}_m$. Ici, $U^+_{\gamma_i}$ est le radical unipotent du
sous-groupe parabolique oppos\'e \`a $P_i$. Posons
$u_i=u U_{\alpha_i,-1}(x)$.

Dans ce cas, $u_i\notin L^{<0}U$. Toutefois, l'unicit\'e de
l'expression $p_i=u_i e_0$ suffit pour l'argument donn\'e dans la
d\'emonstration du lemme 9.1.
Le morphisme qui envoie le point
$$(p_1,\ldots,p_n)\in\prod_{i\notin I(\chi)} (S_{\nu_i}\cap{\bar\Gr}_{\mu_i})
\prod_{i\in I(\chi)} (\phi_{\gamma_i}^{-1}(Uw_iP_i/P_i)\setminus
\epsilon_{\gamma_i}^{-1}(Uw_iP_i/P_i)),$$ sur le point
$$(y_1,\ldots,y_n)\in S_{\nu_\bullet}\cap{\bar\Gr}_{\mu_\bullet}$$ d\'efini par
$$\begin{array}{rcl}
y_1 &=& u_1\vp^{\nu_1} e_0\\
&\ldots & \\
y_n &=& u_1\vp^{\nu_1} u_2\vp^{\nu_2}\ldots u_n\vp^{\nu_n} e_0.\\ \end{array}$$
induit un isomorphisme du sch\'ema source sur un ouvert dense de $C_\chi$.

Pour $i\notin I(\chi)$, on a $u_i\in L^{\geq 0}U$ de sorte que
$\vp^{\sigma_{i-1}} u_i\vp^{-\sigma_{i-1}}$ appartient aussi
\`a $L^{\geq 0}U$ vu que $\sigma_{i-1}$ est dominant.

Pour $i\in I(\chi)$, on a $u_i=u U_{\alpha_i,-1}(x)$ et par cons\'equent,
$$\vp^{\sigma_{i-1}} u_i\vp^{-\sigma_{i-1}} =\vp^{\sigma_{i-1}} u
\vp^{-\sigma_{i-1}}
U_{\alpha_i,\langle\alpha_i,\sigma_{i-1}\rangle-1}(x).$$ Cet \'el\'ement
appartient aussi \`a $L^{\geq 0}U$ parce que
$\langle\alpha_i,\sigma_{i-1}\rangle\geq 1$.

Il s'ensuit que $y_n\in S_\nu\cap {\bar\Gr}_\nu$ de sorte qu'on a un ouvert
dense
de $C_\chi$ contenu dans $m_{\mu_\bullet}^{-1}(S_\nu\cap {\bar\Gr}_\nu)$. Or,
$S_\nu\cap{\bar\Gr}_\nu$ est ferm\'e dans
$S_\nu\cap{\bar\Gr}_{|\mu_\bullet|}$ si bien que la composante $C_\chi$
toute enti\`ere est contenue dans $ m_{\mu_\bullet}^{-1}(S_\nu\cap
{\bar\Gr}_\nu)$. \findem

Il n'est pas difficile de d\'emontrer qu'inversement, si le
$\mu_\bullet$-chemin $\chi$ n'est pas dominant, alors $C_\chi\not\subset
\pi^{-1}(S_\nu\cap\Gr_\nu)$. Nous laissons cette assertion aux soins du
lecteur car elle n'est pas logiquement n\'ecessaire pour la suite. Il nous
suffira seulement de savoir que la multiplicit\'e $\dim(V^\nu_{|\mu_\bullet|})$
de $\calA_\nu$ dans $\calA_{\mu_1}*\cdots*\calA_{\mu_n}$ vaut, au moins, le
nombre de $\mu_\bullet$-chemins dominants allant de $0$ \`a $\nu$.

\begin{proposition}
Pour tout $\lambda\in X^\smv_+$, $\calA_\lambda$ est facteur direct d'un
produit de convolution de la forme
$$
\calA_{\mu_1}*\cdots*\calA_{\mu_n} ,
$$
avec $\mu_1,\ldots,\mu_n\in M$.
\end{proposition}

Compte tenu de 9.4 et de 9.5, il suffit de d\'emontrer qu'il existe un
$\mu_\bullet$-chemin dominant allant de $0$ \`a $\nu$. On d\'emontrera cet
\'enonc\'e combinatoire dans la section 10.

Signalons le corollaire suivant dont on ne se servira pas dans la suite de
l'article. Cet \'enonc\'e se trouve d\'ej\`a dans \cite{Ginzburg} et
\cite{Mirkovic-Vilonen}.

\begin{corollaire}
Pour tous $\lambda,\lambda'\in X^\smv_+$, le produit de convolution
$\calA_{\lambda}*\calA_{\lambda'}$ est un faisceau pervers.
\end{corollaire}

\dem
Si $\calA_\lambda$, resp. $\calA_{\lambda'}$, est un facteur direct de
$\calA_{\mu_1}*\cdots*\calA_{\mu_n}$, resp. de
$\calA_{\mu'_1}*\cdots*\calA_{\mu'_{n'}}$, alors
$\calA_{\lambda}*\calA_{\lambda'}$
est un facteur direct de
$$\calA_{\mu_1}*\cdots*\calA_{\mu_n}*\calA_{\mu'_1}*\cdots*\calA_{\mu'_{n'}} ,$$
qui est pervers d'apr\`es 9.4.
\findem

Cette technique permet aussi de g\'en\'eralier la d\'emonstration de
\cite[cor. 4.3.2]{Ngo2} \`a tout groupe r\'eductif. Cet \'enonc\'e
a \'et\'e d\'emontr\'e auparavant par
Ginzburg, Mirkovic et Vilonen \cite{Ginzburg},
\cite{Mirkovic-Vilonen} du moins lorsque le corps $k$ est le corps des
nombres complexes.

\section{Combinatoire}

Nous proposons deux d\'emonstrations de 9.6.
L'une repose sur un lemme sur les syst\`emes de racines, qui para\^{\i}t
int\'eressant en soi.
L'autre est bas\'ee sur la
th\'eorie des repr\'esentations et le mod\`ele des chemins de Littelmann ;
elle pr\'esente des similitudes remarquables avec des
r\'esultats g\'eom\'etriques des sections 8 et 9, et de cette mani\`ere,
est une bonne illustration de l'\'equivalence tannakienne de
\cite{Ginzburg} et
\cite{Mirkovic-Vilonen}.

\bigskip\noindent{\it Preuve combinatoire.\ }
Rappelons que $M$ d\'esigne l'ensemble des cocaract\`eres minuscules et
quasi-minuscules.
Si $\mu=\gamma^\smv$ est quasi-minuscule,
$\Delta_\gamma$ d\'esigne l'ensemble des racines simples conjugu\'ees \`a
$\gamma$.

Compte tenu de 9.4 et 9.5, la proposition 9.6 d\'ecoule de
l'\'enonc\'e suivant.

\begin{lemme}
Soit $\lambda\in X^\smv_+$. Si $\lambda\not\in M$, il existe une coracine
courte $\beta^\smv$ telle que
$\lambda-\beta^\smv \in X^\smv_+$.
\end{lemme}

\dem Comme $\lambda\not\in M$, il r\'esulte du lemme 1.1
qu'il existe une racine
$\alpha$ telle que $\langle \alpha, \lambda\rangle \geq 2$.
Comme $\lambda$ est dominant, ceci entra\^{\i}ne  que
$\langle \beta, \lambda\rangle \geq 2$, o\`u $\beta$ est la
racine maximale (longue!) du sous-syst\`eme de racines irr\'eductible
de $R$ contenant $\alpha$.
On observe que $\lambda-\beta^\smv$ est un poids de
$V(\lambda)$, donc appartient \`a $\Omega(\lambda)$.

Soit $(\phantom{a}, \phantom{a})$ un produit scalaire $W$-invariant sur
$X^\smv \otimes_{\mathbb Z} {\mathbb R}$,
normalis\'e par la condition que $(\alpha^\smv,\alpha^\smv) = 2$ pour toute
coracine courte $\alpha^\smv$. Pour tout $\chi\in X^\smv$, on posera
$|\chi|^2 := (\chi,\chi)$.

Alors, un calcul facile montre que $\vert \lambda - i\beta^\smv \vert <
\vert \lambda \vert$
pour $0 < i < \langle \beta, \lambda\rangle$. Par cons\'equent, comme
$\langle \beta, \lambda\rangle \geq 2$, on obtient que $\lambda - \beta^\smv$
n'appartient pas \`a $W\lambda$. Notons $\lambda'$ le
cocaract\`ere dominant conjugu\'e \`a
$\lambda-\beta^\smv$, il appartient aussi \`a $\Omega(\lambda)$. On a ainsi
$$
\lambda-\beta^\smv \leq \lambda' < \lambda,
$$
et donc $\lambda' = \lambda - \eta$ et $\beta^\smv = \eta+\nu$, avec
$\eta,\nu\in Q^\smv_+$ et $\eta\not= 0$.

Alors, pour tout $\alpha\in R$, on a
$$
(\lambda,\alpha^\smv)=
\langle\alpha, \lambda\rangle(\alpha^\smv,\alpha^\smv)/2.
$$
et par cons\'equent, $(\lambda,\delta)\geq 0$ pour tout
$\delta\in Q^\smv_+$.

Alors, de l'\'egalit\'e
$$|\lambda-\beta^\smv|^2 = |\lambda-\eta|^2,$$
on d\'eduit que
$$
|\beta^\smv|^2 - |\eta|^2 = 2(\lambda,\nu) \geq 0,
$$
d'o\`u $|\eta|^2 \leq |\beta^\smv|^2 = 2$. Ceci entra\^{\i}ne que $\eta$
est une coracine courte,
et le lemme est d\'emontr\'e. \findem

\bigskip
\noindent{\it Preuve bas\'ee sur la th\'eorie des repr\'esentations.}
Voici un cas tr\`es particulier et bien connu de la r\`egle de
Littlewood-Richard\-son, voir \cite{Littelmann} pour le cas
g\'en\'eral. On rappelle que, pour $\lambda\in X^\smv_+$, $V(\lambda)$
d\'esigne le module simple de plus haut poids $\lambda$
pour le groupe $G^\smv$ d\'efini sur $\Ql$.

\begin{lemme} Soient $\mu\in M$ et $\lambda\in X^\smv_+$.
\begin{enumerate}
\item
 Si $\mu$ est minuscule, alors
$$V(\mu)\otimes V(\lambda)=
\bigoplus_{\deuxind{\nu\in W\mu}{\nu+\lambda\in X^\smv_+}}
V(\lambda+\nu).
$$
\item
Si $\mu$ est quasi-minuscule, alors
$$V(\mu)\otimes V(\lambda)=
\bigoplus_{\deuxind{\nu\in W\mu}{\nu+\lambda\in X^\smv_+}}
V(\lambda+\nu)
\oplus\bigoplus_{\deuxind{\alpha\in\Delta_\gamma}
{\lambda+{\frac 1 2}\alpha^\smv\in X^\smv_+}} V(\lambda).$$
\end{enumerate}
\end{lemme}

\dem
Ceci est bien connu, voir par exemple \cite[Lemma 5A.9]{Jantzen} ou 
\cite[4.2.1]{Donkin} pour le point 1) et  \cite[3.7-3.8]{Polo} pour le 
point 2). 

En utilisant le mod\`ele des chemins de Littelmann \cite{Littelmann}, on
peut aussi argumenter comme suit.
D'apr\`es \loc, le $G^\smv$-module $V(\lambda)$ admet une base
param\'etr\'ee par certains chemins. En particulier, pour $\mu$
minuscule, $V(\mu)$ admet une base
$\{v_{p_{w\mu}} \}_{w\in W/W_\mu}$ o\`u $p_{w\mu}$ est le chemin d\'efini par
$$p_{w\mu}(t)= t\, w\mu,{\rm \ \ pour\ tout\ \ }0\leq t\leq 1;$$
le poids de $v_{p_{w\mu}}$ \'etant $p_{w\mu}(1)=w\mu$.

Pour $\mu=\gamma^\smv$ quasi-minuscule, $V(\mu)$ admet une base
$$\{v_{p_{w\mu}} \}_{w\in W/W_\mu}\cup\{v_{p_\alpha}
\}_{\alpha\in\Delta_\gamma},$$
o\`u pour toute racine simple $\alpha\in\Delta_\gamma$, $p_\alpha$ est le
chemin
$$
p_\alpha(t)= \begin{cases}
-t\alpha^\smv &{\rm \ \ pour\ tout\ \ }0\leq t\leq 1/2 ;\\
(t-1)\alpha^\smv &{\rm \ \ pour\
tout\
\ }1/2\leq t\leq 1;
\end{cases}$$
 le poids de $v_{p_\alpha}$ \'etant $p_\alpha(1)=0$.

D'apr\`es \loc, $V(\mu)\otimes V(\lambda)$ est la somme directe des
$V(\lambda+\chi(1))$, o\`u $\chi$ parcourt l'ensemble des
chemins dans
$V(\mu)$ tels que le translat\'e $\lambda+\chi([0,1])$ soit enti\`erement
contenu dans la chambre dominante. Pour les
chemins $p_{w\mu}$, cela \'equivaut \`a la condition $\lambda+w\mu\in
X^\smv_+$. Pour les chemins $p_\alpha$ avec
$\alpha\in\Delta_\gamma$, cela \'equivaut \`a la condition $\lambda+{\frac
1 2}\alpha^\smv\in X^\smv_+$. \findem

Il r\'esulte de ce lemme que $V(\lambda)$ est un facteur direct d'un
produit tensoriel $V(\mu_1)\otimes\cdots\otimes V(\mu_n)$
si et seulement s'il existe un $\mu_\bullet$-chemin dominant allant de
$0$ \`a $\lambda$. Il suffit donc de d\'emontrer le lemme
suivant.

\begin{lemme}
Pour tout $\lambda\in X^\smv_+$, $V(\lambda)$ est facteur direct d'un
produit tensoriel de la forme
$V(\mu_1)\otimes\cdots\otimes V(\mu_n)$
avec $\mu_1,\ldots\mu_n\in M$.
\end{lemme}

\dem
D\'emontrons d'abord que la repr\'esentation
$$
\rho_M:G\rta\prod_{\mu\in M} \mathop{\rm End}V(\mu)$$ est fid\`ele.
D'abord, il est bien connu, et facile de voir, que pour tout $\xi\in X^\smv_+$,
le sous-groupe de $X^\smv=\Hom(T^\smv,\Gm)$
engendr\'e par les poids de $V(\xi)$
est le sous-groupe engendr\'e par $Q^\smv$ et $\xi$.
D'autre part, on d\'eduit de \cite[Chap.VI, Ex.2.5]{Bourbaki}, que
$M$ contient un syst\`eme de repr\'esentants de $X^\smv/Q^\smv$.
Il en r\'esulte que la restriction de $\rho_M$ au tore maximal $T^\smv$ est
fid\`ele, et donc que $\rho_M$ est fid\`ele.

On en d\'eduit que l'homomorphisme d'alg\`ebres
$${\rm Sym}(\bigoplus_{\mu\in M} V(\mu)\otimes V(\mu)^*)\rta \Ql[G^\smv],$$
o\`u $\Ql[G^\smv]$ d\'esigne
l'alg\`ebre des fonctions r\'eguli\`eres sur $G^\smv$, est surjectif.
D'apr\`es le th\'eor\`eme de Peter-Weyl, tout module $V(\lambda)$
intervient comme facteur direct de l'alg\`ebre $\Ql[G^\smv]$,
d'o\`u le lemme.
\findem

\section{Fin des d\'emonstrations}

On conserve les notations de la section 9. En particulier, soient
$\lambda \in X^\smv_+$, $\nu\in \Omega(\lambda)$ et $\mu_\bullet
= (\mu_1,\dots,\mu_n)$ une suite d'\'el\'ements de $M$ telle que
$\mathcal A_\lambda$ soit un facteur direct de
$\calA_{\mu_1}*\cdots*\calA_{\mu_n}$ (\cf Proposition 9.6).

\bigskip\noindent{\it Preuve du th\'eor\`eme 3.1.}
Compte-tenu des hypoth\`eses ci-dessus,
pour d\'emon\-trer 3.1, il suffit de d\'emontrer que le
complexe $$\rmR\Gamma_c(S_\nu,\calA_{\mu_1}*\cdots*\calA_{\mu_n})$$ est
concentr\'e en degr\'e $2\langle\rho,\nu\rangle$ et que l'endomorphisme de
Frobenius ${\rm Fr}_q$ agit dans $\rmH^{2\langle\rho,\nu\rangle}_c
(S_\nu,\calA_{\mu_1}*\cdots*\calA_{\mu_n})$ comme la multiplication par
$q^{\langle\rho,\nu\rangle}$. D'apr\`es le th\'eor\`eme de changement de
base pour un morphisme propre, on a
$$
\rmR\Gamma_c(S_\nu,\calA_{\mu_1}*\cdots*\calA_{\mu_n})
=\rmR\Gamma_c(m_{\mu_\bullet}^{-1}(S_\nu\cap{\bar\Gr}_{|\mu_\bullet|}),
{\rm IC}({\bar\Gr}_{\mu_\bullet})).
$$
Rappelons qu'on a la stratification
$$
m_{\mu_\bullet}^{-1}(S_\nu\cap{\bar\Gr}_{|\mu_\bullet|})
=\bigcup_{|\nu_\bullet|=\nu}
S_{\nu_\bullet}\cap {\bar\Gr}_{\mu_\bullet}$$
et, d'apr\`es le lemme 9.1, on a un isomorphisme
$$
S_{\nu_\bullet}\cap
{\bar\Gr}_{\mu_\bullet}
\cong(S_{\nu_1}\cap{\bar\Gr}_{\mu_1})\times\cdots\times
(S_{\nu_n}\cap{\bar\Gr}_{\mu_n}),
$$
o\`u $\nu_\bullet = (\nu_1,\dots,\nu_n)$.
De plus, cet isomorphisme est induit par l'isomor\-phisme provenant de
la locale trivialit\'e
$$
\begin{array}{c}
(\vp^{\nu_1}L^{<0}Ge_0\cap {\bar\Gr}_{\mu_1})\times\cdots\times
(\vp^{\nu_n} L^{<0}Ge_0\cap {\bar\Gr}_{\mu_n})\\ \cong
(\vp^{\nu_1}L^{<0}Ge_0\cap
{\bar\Gr}_{\mu_1})\tildetimes\cdots\tildetimes (\vp^{\nu_n}L^{<0}Ge_0\cap
{\bar\Gr}_{\mu_n})
\end{array}
$$
si bien qu'on a
$$
\rmR\Gamma_c(S_{\nu_\bullet}\cap {\bar\Gr}_{\mu_\bullet},
{\rm IC}({\bar\Gr}_{\mu_\bullet}))
=\bigotimes_{i=1}^n
\rmR\Gamma_c(S_{\nu_i}\cap {\bar\Gr}_{\mu_i},\calA_{\mu_i}).
$$
L'assertion \`a d\'emontrer r\'esulte maintenant de 5.2 et de 8.4. \findem

\noindent{\it Preuve du th\'eor\`eme 3.2.} Rappelons que le cas plus facile
$\nu=\lambda$ a \'et\'e d\'emontr\'e dans la discussion qui suit le lemme
5.2. On d\'emontre maintenant le cas plus difficile $\nu\not=\lambda$.

La suite $\mu_\bullet$ a \'et\'e choisie de sorte que la multiplicit\'e
$V^\lambda_{\mu_\bullet}$ de $\calA_\lambda$ dans la d\'ecomposition 9.4 :

$$\calA_{\mu_1}*\cdots*\calA_{\mu_n}=
\bigoplus_{\deuxind{\xi\in X^\smv_+}
{\xi\leq\mu_1+\cdots+\mu_n}}
\calA_\xi\otimes V^\xi_{\mu_\bullet}$$
est non nulle.
On d\'eduit de cette d\'ecomposition l'\'egalit\'e
$$\begin{array}{cl}
& \rmR\Gamma_c(S_\nu,
\calA_{\mu_1}*\cdots*\calA_{\mu_n}\otimes h^*\calL_\psi)\\ = &
\bigoplus\limits_{\deuxind{ \xi\in X^\smv_+}
{\xi\leq\mu_1+\cdots+\mu_n}}
\rmR\Gamma_c(S_\nu,
\calA_\xi\otimes h^*\calL_\psi)\otimes V^\xi_{\mu_\bullet}. \end{array}$$

Du fait que $V^\lambda_{\mu_\bullet}\not=0$ et que $\lambda\not=\nu$, pour
d\'emontrer que
$$\rmR\Gamma_c(S_\nu,
\calA_\lambda\otimes h^*\calL_\psi)=0$$
il suffit de d\'emontrer que la fl\`eche facteur direct $$\begin{array}{cl}
& \rmR\Gamma_c(S_\nu,
\calA_\nu\otimes h^*\calL_\psi)\otimes V^\nu_{\mu_\bullet}\\ \rightarrow &
\rmR\Gamma_c(S_\nu,
\calA_{\mu_1}*\cdots*\calA_{\mu_n}\otimes h^*\calL_\psi) \end{array}
$$
est un quasi-isomorphisme.

Or, d'apr\`es la discussion qui suit le lemme 5.2, on sait que
$$\rmR\Gamma_c(S_\nu,
\calA_\nu\otimes h^*\calL_\psi)\otimes V^\nu_{\mu_\bullet}
=V^\nu_{\mu_\bullet}[-2\langle\rho,\nu\rangle](-\langle\rho,\nu\rangle).$$

Il suffit par cons\'equent de d\'emontrer que pour
$i\not=2\langle\rho,\nu\rangle$, on a $$\rmH^i_c(S_\nu,
\calA_{\mu_1}*\cdots*\calA_{\mu_n}\otimes h^*\calL_\psi)=0$$ et que pour
$i=2\langle\rho,\nu\rangle$, on a
$$\dim(V^\nu_{\mu_\bullet})\geq
\dim(H^{i}_c(S_\nu,
\calA_{\mu_1}*\cdots*\calA_{\mu_n}\otimes h^*\calL_\psi)).$$

Rappelons qu'on a la stratification
$$m_{\mu_\bullet}^{-1}(S_\nu\cap{\bar\Gr}_{|\mu_\bullet|}) =\bigcup_{|\nu_\bullet|=\nu}
S_{\nu_\bullet}\cap {\bar\Gr}_{\mu_\bullet}$$
et que chaque point
$$(y_1,\ldots,y_n)\in S_{\nu_\bullet}\cap {\bar\Gr}_{\mu_\bullet}$$
s'\'ecrit de mani\`ere unique sous la forme $$\begin{array}{rcl}
y_1 & = & x_1\vp^{\nu_1}e_0 \\
&\ldots & \\
y_n & = & x_1\vp^{\nu_1}\ldots x_n\vp^{\nu_n}e_0\\ \end{array}$$
avec $x_i\in L^{<\nu_i}U$ tels que
$x_i\vp^{\nu_i}e_0\in{\bar\Gr}_{\mu_i}$. Pour $\sigma\in X^\smv$, notons
$h_\sigma:LU\rta{\mathbb G}_a$ le morphisme d\'efini par
$h_\sigma(x)=h(\vp^{\sigma}x\vp^{-\sigma})$ ainsi que ses restriction aux
$L^{<\nu}U$ et $S_\nu$. Il est clair que
$$h(y_n)=h(x_1)+h_{\sigma_1}(x_2)+\cdots+h_{\sigma_{n-1}}(x_n).$$ Joint \`a
l'argument de locale trivialit\'e d\'ej\`a utilis\'e dans la preuve de 3.1,
on obtient l'\'egalit\'e $$\begin{array}{cl}
& \rmR\Gamma_c((S_{\nu_\bullet}\cap{\Gr}_{\mu_\bullet}) ,{\rm
IC({\bar\Gr}_{\mu_\bullet})}\otimes h^*\calL_\psi)\\ = & \bigotimes_{i=1}^n
\rmR\Gamma_c((S_{\nu_i}\cap{\Gr}_{\mu_i}) ,{\calA}_{\mu_i}\otimes
h_{\sigma_{i-1}}^*\calL_\psi). \end{array}$$

\begin{lemme}
Si $\sigma\notin X^\smv_+$, alors on a
$${\rm R}\Gamma_c(S_\nu,\calA_\lambda\otimes h_\sigma^*{\mathcal L}_\psi)=0.$$
\end{lemme}

\dem Soit $\alpha\in\Delta$ une racine simple telle que
$\langle\alpha,\sigma\rangle$ soit strictement n\'egatif. Le sous-groupe
${\mathbb G}_a=U_{\alpha,-\langle\alpha,\sigma\rangle-1}$
est alors contenu dans $L^{\geq 0}U$, donc agit de mani\`ere \'equivariante
sur le couple $(S_\nu,\calA_\lambda)$. Or, la restriction de $h_\sigma$ \`a
ce sous-groupe induit l'identit\'e de ${\mathbb G}_a$. Il suffit maintenant
d'appliquer \cite[lemme 3.3]{Ngo}.
\findem

On en d\'eduit l'annulation
$$\rmR\Gamma_c((S_{\nu_\bullet}\cap{\Gr}_{\mu_\bullet}) ,{\rm
IC({\bar\Gr}_{\mu_\bullet})}\otimes h^*\calL_\psi)=0$$ pour les suites
$\nu_\bullet$ dont au moins une des sommes partielles $\sigma_i$ n'est
pas dominante.

Soit maintenant $\nu_\bullet$ une suite avec $\nu_i\in\Omega(\mu_i)$ telle
que toutes les sommes partielles $\sigma_i=\nu_1+\cdots+\nu_i$ sont
dominantes. On dira qu'un $\mu_\bullet$-chemin est de type $\nu_\bullet$
s'il a pour sommets $0,\sigma_1,\ldots,\sigma_n$.
Observons que la condition
$\langle\alpha,\sigma\rangle\geq 1$ apparaissant dans le lemme 8.4 \'equivaut
\`a la condition que $\alpha^\smv/2+\sigma$ soit dominant.

En mettant ensemble les
lemmes 5.3 et 8.4, on arrive \`a l'assertion suivante. Pour
$i\not=2\langle \rho,\nu\rangle $, on a
$$\rmH^i_c(S_{\nu_\bullet}\cap{\bar\Gr}_{\mu_\bullet}, {\rm
IC}({\bar\Gr}_{\mu_\bullet})\otimes h^*\calL_\psi)=0$$
et pour $i=2\langle \rho,\nu\rangle $, on a
$$\begin{array}{cl}
& \dim(\rmH^i_c(S_{\nu_\bullet}\cap{\bar\Gr}_{\mu_\bullet}, {\rm
IC}({\bar\Gr}_{\mu_\bullet})\otimes h^*\calL_\psi))\\
= & |\{\mu_\bullet\mbox{-chemins dominants de type\ }\nu_\bullet\}|.
\end{array}$$

Par ailleurs, compte tenu de 9.4 et de 9.5, on a l'in\'egalit\'e
$$\dim(V^\nu_{\mu_\bullet})\geq
|\{\mu_\bullet \mbox{-chemins dominants allant de 0 \`a\ } \nu\}|.$$
La d\'emonstration du th\'eor\`eme 3.2 est termin\'ee. \findem

\renewcommand{\refname}{R\'ef\'erences}

\begin{small}

\end{small}
\bigskip
\bigskip
{CNRS, UMR 7539, LAGA,\\
Institut Galil\'ee, \\
Universit\'e Paris-Nord,\\
93430  Villetaneuse, France.

\bigskip\noindent{\it Courriers \'electroniques :}
\begin{verbatim}
ngo@math.univ-paris13.fr
polo@math.univ-paris13.fr
\end{verbatim}}

\end{document}